\documentclass[3p,times,11pt]{elsarticle}

\usepackage{amssymb}

\usepackage{lineno}

\usepackage{amssymb,amsmath,mathrsfs}
\usepackage{graphicx,epstopdf,subcaption}
\usepackage[dvipsnames]{xcolor}
\usepackage{url}		
\usepackage{algorithm} 
\usepackage{algorithmicx} 
\usepackage{algpseudocode}
\usepackage{float}
\usepackage{tikz}
\usetikzlibrary{patterns}
\usepackage{mathtools}
\mathtoolsset{showonlyrefs=true}
\usepackage{multirow}

\newcommand{\dee}{ {\,\,\mbox{d}} }

\newcommand{\bfE}{ \mathbf{E} }

\newcommand{\bfA}{\mathbf{A}}
\newcommand{\bfH}{\mathbf{H}}
\newcommand{\bfT}{\mathbf{T}}
\newcommand{\bfI}{\mathbf{I}}

\newcommand{\calB}{\mathcal{B}}

\newcommand{\calF}{\mathcal{F}}

\newcommand{\calL}{\mathcal{L}}

\newcommand{\bff}{ \mathbf{f} }

\newcommand{\bfn}{ \mathbf{n} }
\newcommand{\bfs}{ \mathbf{s} }

\newcommand{\bfu}{ \mathbf{u} }

\newcommand{\bfv}{ \mathbf{v} }

\newcommand{\bfx}{ \mathbf{x} }

\newcommand{\bbE}{ \mathbb{E} }

\newcommand{\bbP}{ \mathbb{P} }
\newcommand{\bbR}{ \mathbb{R} }

\newcommand{\etalow}{\eta_{\mathrm{low}}}
\newcommand{\etaupp}{\eta_{\mathrm{upp}}}

\newcommand{\etamlow}{\eta_{\mathrm{low}}^{-}}
\newcommand{\etamupp}{\eta_{\mathrm{upp}}^{-}}

\newcommand{\etaplow}{\eta_{\mathrm{low}}^{+}}
\newcommand{\etapupp}{\eta_{\mathrm{upp}}^{+}}

\newcommand{\etapmlow}{\eta_{\mathrm{low}}^{\pm}}
\newcommand{\etapmupp}{\eta_{\mathrm{upp}}^{\pm}}

\newcommand{\zetaupp}{\zeta_{\mathrm{upp}}}

\newcommand{\xiupp}{\xi_{\mathrm{upp}}}

\newcommand{\ds}{\displaystyle}

\newcommand{\red}[1]{\textcolor{red}{#1}}

\newtheorem{thm}{Theorem}

\newtheorem{rmk}{Remark}

\journal{ }

\newcommand{\rev}[1]{#1}
\newcommand{\prev}[1]{{\color{black}{#1}}}

\begin{document}


\begin{frontmatter}




\title{Goal-Oriented Adaptive Modeling of Random Heterogeneous Media and Model-Based Multilevel Monte Carlo Methods
}


\author{Laura Scarabosio$^1$, 
Barbara Wohlmuth$^1$,
J. Tinsley Oden$^2$, Danial Faghihi$^2$}

\address{
$^1$ Technical University of Munich, Germany, Department of Mathematics, Chair of Numerical Mathematics (M2)\\
$^2$ Institute for Computational Engineering and Sciences (ICES), The University of Texas at Austin
}


\begin{abstract}

Methods for generating sequences of surrogates approximating fine scale models of two-phase random heterogeneous media are presented that are designed to adaptively control the modeling error in key quantities of interest (QoIs). For specificity, the base models considered involve stochastic partial differential equations characterizing, for example, steady-state heat conduction in random heterogeneous materials and stochastic elastostatics problems in \prev{linear} elasticity. The adaptive process involves generating a sequence of surrogate models defined on a partition of the solution domain into regular subdomains and then, based on estimates of the error in the QoIs, assigning homogenized effective material properties to some subdomains and full random fine scale properties to others, to control the error so as to meet a preset tolerance.
New model-based Multilevel Monte Carlo (mbMLMC) methods are presented that exploit the adaptive sequencing and are designed to reduce variances and thereby accelerate convergence of Monte Carlo sampling. Estimates of cost and mean squared error of the method are presented.
The results of several numerical experiments are discussed that confirm that \prev{substantial saving} in computer costs can be realized through the use of controlled surrogate models and the associated mbMLMC algorithms.
\end{abstract}


\end{frontmatter}

\section{Introduction}

The use of high-dimensional, high-fidelity computational models to simulate complex physical phenomena in heterogeneous material systems has increased in recent years. This has been due, in part, to substantial increases in computational power that occurred during this period, but also due to the increased reliance on computer predictions as a basis for critical decisions effecting the design and performance of complex engineered systems. 
It is now more frequently recognized that predictive models of such system must account for uncertainties in microstructural properties and that stochastic models of random media are often required for reliable predictions. 
Paradoxically, models which can deliver the desired resolution of the critical quantities of interest (the QoIs), while also capturing the stochastic character, are often of such size and complexity that the desired solutions are  intractable. So, methods for reducing the cost of computations by reducing the model size while also retaining accuracy in the prediction of QoIs are of great value.

In this paper, we address this broad problem by developing methods to generate sequences of tractable stochastic surrogate models less complex and of smaller dimension than \rev{fine scale} models, that deliver accurate approximations of key \prev{QoIs for random, heterogeneous media. We derive estimates of modeling error in specific output functionals and use these to control adaptively the choice of our surrogates. Our adaptive modeling process provides a framework for developing new Monte Carlo \prev{(MC)} based solution paradigms under the broad category of Multilevel Monte Carlo \prev{(MLMC)} methods. We consider two classes of problems: diffusion in a random two-phase medium and elastostatics of composite materials with randomly distributed isotropic inclusions.}

The idea of goal-oriented estimation and control of modeling error was introduced \prev{by Oden and Vemaganti in connection with multiscale modeling of heterogeneous materials in \cite{oden1999adaptive, odenvemaganti2000,vemagantioden2001}, and generalized to nonlinear problems by Oden and Prudhomme \cite{odenprudhomme2002}}. The reference to \textit{goal-oriented estimates} is intended to imply local estimates of \prev{QoIs}, generally characterized by linear functionals of the solution to the forward problem, as opposed to global estimates of errors in various norms.
Modeling error, as opposed to discretization error, is the relative error in the solutions or \prev{QoIs} between low-fidelity or coarse-scale approximations and those of a high-fidelity base model, which, while generally intractable, is assumed to provide an  high-resolution characterization of the phenomenon of interest.
In this exposition, we focus on the estimation and control of modeling error and assume that discretization errors are made negligible by using sufficiently refined meshes. A general framework for estimating modeling error is described in \cite{odenprudhomme2002,prudhommeoden1999} and methods for adaptive modeling, including hybrid continuum-atomistic models, molecular dynamics models, and models of nonlinear systems, are discussed in \cite{odenprudrombaum2006,baumanodenprud2009}. Further generalizations and applications of these methods to broad classes of problems were given in Oden et al \cite{odenprudrombaum2006}.
Braack and Ern \cite{braack2003} described techniques for a posteriori estimation and control of both modeling and discretization error. \prev{The construction of surrogate models with error control within a stochastic framework is also described in \cite{Chamoinetal} for atomic-to-continuum modeling and \cite{Zaccardi} for a stochastic-deterministic coupling method.} More recently, Maier and Rannacher \cite{maier2016arx,maier2016} continued a duality-based approach to model adaptivity for heterogeneous materials that involves a post-processing procedure for selecting optimal models for estimating errors in \prev{QoIs}.

In the present paper, we describe extensions of the notions of modeling error estimation and control in which the underlying stochastic high-fidelity model is characterized by stochastic elliptic boundary value problems. A large literature exists on the physics of such random heterogeneous media; see for example, the treatises of Torquato  \cite{torquato2013} and Buryanchenko \cite{buryachenko2007}, and the survey of Jeulin and Ostoja-Starzewski \cite{jeulinostoja2001}. 
The central aim of most theories of random media is to derive effective properties of homogenized representation of the material \cite{jikovkozlov2012,
palencia1983,
bensoussan1978}, which is also a factor in the adaptive algorithms developed here.

\prev{The sequence of models that we construct is employed in a MLMC setting. The classical MLMC method \cite{giles2008, giles_2015} uses a hierarchy of space (or time) discretizations of a partial (resp. ordinary) differential equation to accelerate the Monte Carlo convergence on the finest level. 
Its efficiency relies on a delicate trade-off between computational cost and variance reduction across the levels, which is often guaranteed by convergence theorems for the discretized solution. This is not the case for a sequence of surrogate models, where the unavailability of convergence rates (at least a priori) constitutes the main challenge in the construction of an efficient multilevel estimator. 
The use of sets of low-fidelity models correlated with a high-fidelity base model as a tool for deriving \prev{MC} solvers for large-scale simulations has led to the development of Multi-Fidelity Monte Carlo Methods (MFMC) \cite{willcoxgunzburger2016I,willcoxgunzburger2016II,willcox2017}.  In these methods, the selection of surrogate models and the optimal distribution of samples among the levels are tied to the correlation between the approximate QoIs of each surrogate and the QoI of the high-fidelity model, and by the cost of evaluating the QoI for each surrogate model. 
The procedure proposed in \cite{willcoxgunzburger2016I} is independent of the availability of a priori or a posteriori error estimators. Instead, our goal is to design a MLMC algorithm that exploits the information provided by our a posteriori error estimator. 
We show that the latter allows the ordering of models with respect to their accuracy and enables one to fit this model hierarchy into the framework of standard MLMC. 
A similar strategy has been adopted in \cite{vidal2016, vidal2015}, but in these works the surrogate models correspond to reduced basis approximations of the high-fidelity model and no a posteriori error estimator is needed for constructing and ordering them. We remark that, in the case of a hierarchy based on mesh refinement, the MFMC estimator of \cite{willcoxgunzburger2016I} has been shown to provide the same performance as the MLMC estimator, while distributing differently the samples across the levels \cite{peherstorfer2016cvg}. 
A combination of low and high-fidelity models has also been explored in the stochastic collocation framework in \cite{narayan2014mfsc,zhu2017mfmc}.}

Our approach addresses several challenging problems, including
1) the derivation of \textit{two-sided bounds} on goal-oriented estimates of modeling error in QoIs in systems modeled by stochastic PDEs;
2) the construction of homogenized models of the microstructure of random two-phase media;
3) the development and implementation of an adaptive modeling algorithm to control the error in QoIs delivered by models representing a mixture of homogenized and microscale media; and
4) the solution of representative model problems to demonstrate the implementation and effectivity of the proposed methodology.
\prev{Related works that address these questions are \cite{Romkesetal} and \cite{Zaccardi}, but there surrogate models are constructed by coupling a deterministic homogenized model to a local stochastic model, and statistics are computed by Monte Carlo integration. The main novelty of our work is to leverage an adaptive strategy to contruct surrogate models corresponding to levels in a MLMC algorithm. We are able to demonstrate improvements in the efficiency of the MLMC algorithm compared to plain \prev{MC} when estimates of modeling errors are used to guide MC-type solvers. In order to achieve variance reduction across the levels, differently from \cite{Romkesetal, Zaccardi} we use surrogates which are stochastic in the whole computational domain.}

In Section 2, following this introduction, the notation and mathematical structure of a class of stochastic boundary-value problems are presented.
The idea of constructing sequences of lower-dimensional surrogate models approximating the high-fidelity base model is taken up in Section \ref{sec:surrpb}.
A posterior estimates of error in QoIs are presented in Section \ref{sec:goalerrest}.
Applications of the methodologies to \prev{heat transfer and linear elasticity, as well as} the adaptive modeling algorithm for controlling the error are discussed in Section \ref{sec:heat}. 
The theory and algorithms underlying MLMC methods based on levels generated by sequences of surrogates, which we refer to as model-based MLMC (or mbMLMC), are presented in Section \ref{sec6:MLMC}.
Numerical experiments involving applications of the theory and methodology to model problems are given in Section \ref{sec:numexp}.
It is demonstrated that the use of a posterior estimates of error in \prev{QoIs} delivered by a sequence of surrogate models can accelerate the convergence of MC methods. 
Concluding comments are collected in Section \ref{sec:conclusions}.

\section{Stochastic Models Involving Random Heterogeneous Media}

We consider a class of linear stochastic elliptic boundary-value problems that model various physical phenomena in random, heterogeneous media. We assume the media occupies a bounded domain $D$ in $\bbR^d, d = 1, 2, 3$, with Lipschitz boundary $\partial D$. The problem is set in a complete probability space $(\Omega,\calF,\bbP)$, \prev{where} the sample set $\Omega$ of possible outcomes  describes realizations of microstructural distributions of the two-phase material, 
$\calF$ the $\sigma$-algebra of subsets of $\Omega$, and
$\bbP$ the probability measure.
The problem is to find a stochastic field 
$\mathbf{u}: \Omega \times \bar{D} \rightarrow \mathbb{R}^m, m=1,2,3$ and compute a QoI from it.  \prev{For ease of presentation, we set $m=1$ and denote by $u=u(\omega,\bfx)$, $\omega\in\Omega$, $\bfx\in D$, the scalar stochastic field. We} first consider a model class of problems formally governed by the equations
\begin{linenomath}
\begin{equation}
\arraycolsep=1.4pt\def\arraystretch{1.4}
\left.\begin{array}{rcll}
\nabla \cdot \bfA(\omega,\bfx)\nabla u(\omega,\bfx) & = & f(\omega,\bfx), & \bfx \in D,\\
\bfA(\omega,\bfx)\nabla u(\omega,\bfx) \cdot\bfn & = & \sigma(\omega,\bfx), & \bfx\in\Gamma_N,\\
 u(\omega,\bfx) & = & \mathbf{0}, & \bfx\in\Gamma_D,
\end{array}\right\}
\label{eq:strongstochastic}
\end{equation}
\end{linenomath}
almost surely for all $\omega \in \Omega$, where $\boldsymbol{\nabla}$ is the spatial gradient operator,
$\mathbf{A}(\omega,\mathbf{x})$ is a $d \times d$ symmetric positive definite matrix representing the conductivity, diffusivity, mobility, etc. of the random field,
$f(\omega, \bfx)$ and $\sigma(\omega, \bfx)$ are prescribed data,
$\Gamma_N$ being a subset of $\partial D$ on which Neumann data is prescribed, and
$\Gamma_D = \partial D\setminus \Gamma_N$. Later we focus on the special case $\bfA(\omega, \bfx) = \kappa(\omega,\bfx)\bfI$. \prev{The presentation in this and the next sections generalizes naturally to the elasticity case, that is also considered, and for which $\bfA(\omega, \bfx) = \bfE(\omega,\bfx)$ is the fourth order elasticity tensor and $\bfu = \bfu(\omega,\bfx)$ the vector-valued stochastic displacement field.}

If $\calB(D)$ denotes the Borel $\sigma$-algebra generated by open sets of $D$,
then $\bfA$, $f$ and $\sigma$, as $(L^{\infty}(D))^{d\times d}$-, \prev{$H^{-1}_{\Gamma_D}(D)$-,}
and $H^{-\frac{1}{2}}(\Gamma_N)$-valued quantities, respectively, are assumed measurable with the induced $\sigma$-algebra, $\calF\otimes\calB (D)$ for
$\bfA$ and $f$, and $\calF \otimes\calB(\Gamma_N)$ for $\sigma$, $\calB(\Gamma_N)$
being the sigma-algebra
associated with open sets in $\Gamma_N$.  \prev{Here, $H^{-1}_{\Gamma_D}(D)$ denotes the dual of  $H^{1}_{\Gamma_D}(D)= \{ v \in H^1(D), {\rm tr}v|_{\Gamma_D}=0 \}$.}

For almost every $\omega \in \Omega$, $\mathbf{A}(\omega, \bfx)$ is assumed to be bounded and coercive, \prev{uniformly} with respect to $\mathbf{x}$ and $\omega$,
$\mathbf{A}(\omega, \bfx)$ representing a possibly \rev{discontinuous, highly oscillatory} function characterizing the two phases of the media.
Thus, there exist $\rev{0<}\alpha_{min},\alpha_{max}\rev{<+\infty}$ such that
\begin{linenomath}
\begin{equation}
\bbP(\omega\in\Omega : \alpha_{min}\mathbf{a}^T\mathbf{a} \le
\mathbf{a}^T\mathbf{A}(\omega, \bfx)\mathbf{a} \le
\alpha_{max}\mathbf{a}^T\mathbf{a},
\forall \mathbf{a}\in\bbR^d, \mathrm{ a.e. } \; \bfx\in \bar{D}) =1.
\label{eq:eminmax}
\end{equation}
\end{linenomath}

\subsection{Function Spaces}
For an $\bbR^N$-valued random function $Y \in L^1_{\bbP}(\Omega)$, its expectation
is denoted by 
\begin{equation}
\bbE[Y] := \int_{\Omega} Y (\omega) \dee\bbP(\omega) = 
\int_{\bbR^N} y \dee \mu_{Y} (y),
\end{equation}
where $\mu_{Y}$ is the distribution measure for $Y$. We assume $\mu_{Y}$
is absolutely continuous with respect to the Lebesgue measure, so there exists
a probability density function $\pi:\bbR \rightarrow[0,\infty)$ such that  
\begin{equation}
\bbE[Y] = \int_{\bbR^N} y \,\, \pi(y) \dee y.
\end{equation}

\prev{For a Hilbert space $V$, we denote by $L^2_{\bbP}(\Omega,V)$ the space of square-integrable, $V$-valued stochastic functions, equipped with the norm $\|v\|_{L^2_{\bbP}(\Omega,V)} := \left\{\bbE\left[\lVert v\rVert_V^2\right]. 
\right\}^{\frac{1}{2}}$. It holds $L^2_{\bbP}(\Omega,V)\simeq L^2_{\bbP}(\Omega) \otimes V$, cf.~\cite{babuskatempone2004}.}

\subsection{Variational Form}
The weak or variational formulation of the stochastic diffusion problem
\eqref{eq:strongstochastic} is defined as follows: let $H :=
L^2_{\bbP}(\Omega)\otimes H^1_{\Gamma_D}(D)$ denote the tensor product Hilbert
space endowed with the inner product
\begin{linenomath}
\begin{equation*}
(u,v)_H := \bbE\left[\ds\int_D (\nabla u\cdot\nabla v)\dee \bfx\right] = \ds\int_{\Omega}\ds\int_D \nabla u(\omega,\bfx) \cdot \nabla v(\omega,\bfx) \dee \bfx \dee \bbP(\omega)
\nonumber.
\end{equation*}
\end{linenomath}
We define the bilinear form and linear form,
\begin{linenomath}
\begin{equation}
B(u,v) = \bbE\left[  \int_D \bfA\nabla u \cdot \nabla v\dee\bfx \right],
\,\,\,\,\,B: H\times H\rightarrow \bbR,
\label{eq:Bdefn}
\end{equation}
\begin{equation}
F(v) = \bbE\left[ \int_D  f\cdot v\dee\bfx + \int_{\Gamma_N}
 \sigma\cdot v\dee\bfs \right],
\,\,\,\,\,F: H\rightarrow \bbR.
\label{eq:Fdefn}
\end{equation}
\end{linenomath}
By virtue of \eqref{eq:eminmax}, $B(\cdot,\cdot)$ is continuous and coercive.
We assume that $f\in L^2(\Omega,L^2(D))$ and $\sigma \in
L^2(\Omega,H^{-\frac{1}{2}}(\Gamma_N))$, so that $F(\cdot)$ is continuous.

The stochastic variational formulation of \eqref{eq:strongstochastic} is then:
\begin{linenomath}
\begin{equation}
\left.\begin{matrix}
\mbox{Find } u \in H \mbox{ such that}\\[0.5em]
B(u,v) = F(v),\,\,\mbox{ for all } v \in H. \\ 
\end{matrix}\right\}
\label{eq:weakstochastic}
\end{equation}
\end{linenomath}
Under the assumptions laid down so far, a unique solution to
\eqref{eq:weakstochastic} exists by the Lax-Milgram theorem.

Upon solving \eqref{eq:weakstochastic} for the stochastic field $u$, we 
wish to evaluate specific {\em quantities of interest} (QoIs), the goals of
the modeling and simulation, represented by continuous linear functionals 
\begin{linenomath}
\begin{equation}
Q(u) = \mathbb{E}[q(u)] \in \mathbb{R},\quad q: H^1_{\prev{\Gamma_D}}(D) \; \rightarrow \mathbb{R}.
\end{equation}
\end{linenomath}
Many examples of meaningful QoIs could be cited and typical cases are given later. It should be emphasized
that \textit{the calculation of QoIs is the principal goal of constructing the
mathematical model \eqref{eq:strongstochastic} and of solving it numerically.}
The actual solution $u(\omega,\bf x)$ is of interest only as a step in
computing the target QoIs. This fact is central in our approach to the
analysis of random heterogeneous media. \rev{While we addess the forward problem, goal-oriented approaches for Bayesian inverse problems can be found, for instance, in \cite{mattis2018goal} for parameter and mesh adaptivity and in \cite{prudhomme2015adaptive} for model adaptivity in turbulence simulations.}

Corresponding to each QoI, $Q\in\calL(H,\bbR)$, is a unique function
$w = w(\omega,\bfx) \in H$ that represents the stochastic
field generated by the input $Q$, and which is defined as the solution of the following adjoint or
dual problem,
\begin{linenomath}
\begin{equation}
\left.\begin{matrix}
\mbox{Given } Q:H\rightarrow\bbR, \mbox{ find } w \in H \mbox{ such that}\\[0.5em]
B(v,w) = Q(v),\,\,\mbox{ for all } v \in H.
\end{matrix}\right\}
\label{eq:weakadjoint}
\end{equation}
\end{linenomath}
The solution $w$ can be interpreted as the generalized
Green's function corresponding to the functional $Q$.

\section{Surrogate Problem Classes}\label{sec:surrpb}

Unfortunately, \rev{problem
\eqref{eq:weakstochastic} (or \eqref{eq:strongstochastic}) and problem
\eqref{eq:weakadjoint} are generally
hardly tractable in real-life applications, owing to their enormous size and complexity.} 
An alternative approach must be
explored in which \eqref{eq:weakstochastic} (or \eqref{eq:strongstochastic})
and \eqref{eq:weakadjoint} are replaced by a tractable sequence of
approximations that can be solved \rev{much cheaply} and which produce sufficiently accurate 
approximations of the QoIs. 

Therefore, we introduce the surrogate primal and adjoint problems
\begin{linenomath}
\begin{equation}\label{eq:surg_forward_adj}
\left.\begin{matrix}
\mbox{Find } u_0 \in H, \mbox{ such that}\\[0.5em]
B_0(u_0,v) = F(v) \; {\rm for \; all} \;\; v \in H,
\end{matrix}\right\}\qquad\quad
\left.\begin{matrix}
\mbox{Find } w_0 \in H, \mbox{ such that}\\[0.5em]
B_0(v,w_0) = Q(v) \;{\rm for \; all} \;\; v \in H,
\end{matrix}\right\}
\end{equation}
\end{linenomath}
where 
\begin{linenomath}
\begin{equation}
B_0(u_0,v) = \bbE\left[  \int_D \bfA_0\nabla u_0 \cdot \nabla v
\dee\bfx \right],
\end{equation}
\end{linenomath}
and $\bfA_0$ is a suitable approximation of $\bfA(\omega,\bfx)$.
In Section \ref{sec:heat}, we describe how to derive sequences of surrogate pairs $\{(u_0)^{(l)}), (w_0)^{(l)}) \}$, $l= 1, 2, \cdots, L$, of solutions to \eqref{eq:surg_forward_adj} designed to adaptively yield approximations of the QoI, $Q(u)$, of increasing accuracy.
Generally, $\bfA_0(\omega,\bfx)$ will be of the form,
\begin{linenomath}
\begin{equation}
\bfA_0(\omega,\bfx) = \begin{cases}
\bfA(\omega,\bfx), &\bfx \in D_{fine},\\
\bfA_{00}(\omega,\bfx), &\bfx \in D\setminus D_{fine},
\end{cases}
\label{eq:E0kdefn}
\end{equation}
\end{linenomath}
where $\bfA(\omega,\bfx)$ is the ``fine scale'' or heterogeneous coefficient matrix in \eqref{eq:strongstochastic},
$D_{fine}$ is a subset of $D$ containing fine scale random features of the microstructure, and $\bfA_{00}(\prev{\omega},\bfx)$ characterizes homogenized material features.

\section{Goal-Oriented A Posteriori Estimates of Modeling Error}\label{sec:goalerrest}

In this section, we present an extension of the analysis in \cite{odenvemaganti2000,vemagantioden2001} to the stochastic
systems \eqref{eq:strongstochastic} (or \eqref{eq:weakstochastic} and \eqref{eq:weakadjoint}). We note that the positive definite bilinear form $B(\cdot,\cdot)$ of \eqref{eq:Bdefn}
generates an inner product and an energy norm \prev{on $H$, given by $\|v\|_B := \sqrt{B(v,v)}$,}
which, by virtue of \eqref{eq:eminmax}, is equivalent to the norm
$\|v\|_{H} = \sqrt{(v,v)_H}$, with equivalence constants independent
of $\omega\in\Omega$.
Also, as noted earlier, we consider only modeling error under the assumptions that discretization errors are negligible.

The following result, \prev{closely related to \cite[Thm. 6]{Romkesetal},} establishes computable two-sided error bounds on the \prev{QoI}. We make use of the notation $\bfI_0 = \bfI_{0}(\omega,\mathbf{x}):= \bfI - \bfA^{-1}(\omega,\mathbf{x})\bfA_0(\omega,\mathbf{x})$.

\begin{thm}\label{thm:errest}
Let $u_0$ and $w_0$ be solutions to the surrogate primal and adjoint problems
\eqref{eq:surg_forward_adj}, respectively. Then the modeling error in the
quantity of interest is bounded above and below as follows:
\begin{linenomath}
\begin{equation}\label{eq:bound}
\etalow \le Q(e_0) \le \etaupp,
\end{equation}
\end{linenomath}
where
\begin{linenomath}
\begin{equation}
\etalow := \frac{1}{4} (\etaplow)^2 - \frac{1}{4}(\etamupp)^2 + 
{R}_{u_0}(w_0),
\label{eq:lofeteqnlhs}
\end{equation}
\begin{equation}
\etaupp := \frac{1}{4} (\etapupp)^2 - \frac{1}{4}(\etamlow)^2 +
{R}_{u_0}(w_0),
\label{eq:lofeteqnrhs}
\end{equation}
\end{linenomath}
with
\begin{linenomath}
\begin{equation}
\etapmlow := \frac{| {R}_{su_0\pm s^{-1}w_0}( u_0 + 
\theta^{\pm} w_0)|}{\| u_0 + \theta^{\pm} w_0\|_{{B}}},
\label{eq:etatpmldef}
\end{equation}
\end{linenomath}
$s \in \mathbb{R}$ being real numbers, 
and
\begin{linenomath}
\begin{equation}
(\etapmupp)^2 := s^2 \zetaupp^2
+ s^{-2}\xiupp^2
\pm \bbE\left[
2 \int_{D} {\bfI_0} \nabla u_0\cdot \bfA\,{\bfI_0} \nabla
w_0 \dee \bfx \right],
\label{eq:etatpmudef}
\end{equation}
\begin{equation}
\theta^{\pm} = 
\frac{B(u_0,w_0) R_{u0}(su_0 \pm s^{-1}w_0) - B(u_0,u_0) R_{w_0}(su_0 \pm s^{-1}w_0)}
{B(u_0,w_0) R_{w0}(su_0 \pm s^{-1}w_0) - B(w_0,w_0) R_{u_0}(su_0 \pm s^{-1}w_0)},
\end{equation}
\begin{equation}
\zetaupp = \| u_0 \|_{\mathbf{I}_0B}, \qquad
\xiupp = \| w_0 \|_{\mathbf{I}_0B},
\end{equation}
$\| \cdot \|_{\mathbf{I}_0B}$ being the $\mathbf{I}_0$-weighted norm,
\begin{equation}
\| u_0 \|_{\mathbf{I}_0B} = 
\bbE\left[\ds\int_{D} \bfA\,{\bfI_0} \nabla u_0\cdot \nabla u_0
\dee \bfx \right],
\end{equation}
\end{linenomath}
and $R_{u_0}(\cdot)$ is the linear functional, defined, for any $g\in H$, by
\begin{linenomath}
\begin{equation}\label{eq:rgv}
R_g(v) = 
\bbE\left[\ds\int_{D} \bfA\,{\bfI_0} \nabla g\cdot \nabla v
\dee \bfx \right].
\end{equation}
\end{linenomath}
\qed
\end{thm}

The proof follows from a straightforward generalization of the \prev{proofs in \cite{oden1999adaptive,Romkesetal}.}

It has been noted in the numerical experiments of \cite{odenvemaganti2000} that $\eta_{upp}$ and $\eta_{low}$ can provide poor estimates of the error. Simpler computable but \rev{effective} approximate error estimators can be derived from the estimators in \eqref{eq:lofeteqnlhs} and \eqref{eq:lofeteqnrhs} and used to guide the adaptive processes described in the next sections. 

Considering the estimate of the modeling error in the QoI based only on the upper bounds, we obtain the approximate estimator
\begin{linenomath}
\begin{equation}
\eta_{est} :=  \frac{1}{4} (\eta^+_{upp})^2 - \frac{1}{4} (\eta^-_{upp})^2 +  R_{u_0}(w_0)= \mathbb{E}\left[-\int_D \bfA^{-1}\bfA_0 \bfA\bfI_0\nabla u_0\cdot\nabla w_0 \dee\bfx\right]. \label{eq:eta_est}
\end{equation}
\end{linenomath}

A similar second error estimator is obtained as $\eta_{est,low}:=\frac{1}{4} (\eta^+_{low})^2 - \frac{1}{4} (\eta^-_{low})^2 +  R_{u_0}(w_0)$. 

\section{Heat Transfer \prev{and Linear Elasticity} in Two-Phase Materials with Random Microstructure}\label{sec:heat}

In this section, we focus on the case of a two-phase isotropic random media. \prev{In Subsection \ref{ssec:heat} we describe the construction of surrogate models and the a posteriori error estimator for the heat transfer problem, while in Subsection \ref{ssec:ela} we present the extension to linear elasticity}.

\subsection{Adaptive Modeling and Construction of Surrogate Models \prev{for Heat Transfer}}\label{ssec:heat}
\prev{In this case,} $\mathbf{A}$ is of the form $\mathbf{A}(\omega,\mathbf{x}) = \kappa(\omega,\mathbf{x})\mathbf{I}$, 
$\mathbf{I}$ being the identity, and $\kappa$ the \prev{random field describing material conductivities varying between two phases. Denoting by $\kappa_M$ the conductivity for the matrix material and by $\kappa_I$ the one for the inclusions, formally
\begin{linenomath}
\begin{equation}\label{eq:kappaIM}
\kappa (\omega,\mathbf{x}) = \kappa_M\mathcal{X}_M(\omega,\mathbf{x}) + \kappa_I\mathcal{X}_I(\omega,\mathbf{x}), \qquad \omega\in\Omega, \quad \mathbf{x}\in\bar{D},
\end{equation}
\end{linenomath}
$\mathcal{X}_M$ and $\mathcal{X}_I$ being the characteristic functions for the two phases. Then the solution $u(\omega,\mathbf{x})$ to \eqref{eq:weakstochastic} is the corresponding temperature field.}

\prev{The} bilinear and linear forms of \eqref{eq:Bdefn} and \eqref{eq:Fdefn} reduce to,
\begin{linenomath}
\begin{equation}\label{eq:bilinheat}
B(u,v) = \bbE\left[  \int_D \kappa(\cdot,\bfx) \nabla u(\cdot,\bfx) \cdot \nabla v(\cdot,\bfx)\dee\bfx \right],
\end{equation}
\begin{equation}\label{eq:linheat}
F(v) = \bbE\left[ \int_D  f(\cdot,\bfx)\cdot \dee\bfx + \int_{\Gamma_N}
 \sigma(\cdot,\bfx)\cdot \mathbf{n}\dee\bfx \right].
\end{equation}
\end{linenomath}

The QoI is always assumed to be a local feature of the physical system fully characterized by the solution of the \rev{fine scale} forward problem. Examples include the average temperature over a subdomain $A_q\subset D$,
\begin{linenomath}
\begin{equation}\label{eq:sur_qoi1}
Q(u(\omega,\mathbf{x})) = \mathbb{E} \left[ \frac{1}{| A_q |} \int_{A_q} u(\omega,\mathbf{x}) d\mathbf{x} \right],
\end{equation}
\end{linenomath}
$| A_q |$ being the \prev{volume} of $A_q$,
or \rev{the heat flux through the boundary of a subdomain.}
%

Following the ideas in \cite{odenvemaganti2000}, we partition the domain $D$ into blocks $\left\{B_1,\ldots,B_{K_b}\right\}$, $K_b\in\mathbb{N}$, so that $D=\bigcup_{k=1}^{K_b}\overline{B}_k$; see for instance Figure \ref{fig:Lshapesmooth} in Section 8. As an approximation of the conductivity $\kappa$, we consider
\begin{linenomath}
\begin{equation}
\kappa_0(\omega,\bfx) = 
\begin{cases}
\kappa(\omega,\bfx)&\rev{ \text{if } \bfx\in B_k\text{ such that }B_k\subseteq D_{fine},}\\
\kappa_{eff} (\omega,\rev{k})&\rev{ \text{if }\bfx \in B_k\text{ such that }B_k\subseteq D\setminus D_{fine},}
\end{cases}
\label{eq:kappa0}
\end{equation}
\end{linenomath}
where $D_{fine}$ consists of a collection of blocks around the QoI and $\kappa_{eff}(\omega,\cdot)$, $\omega\in\Omega$, is a \textsl{blockwise}-homogenized coefficient. \prev{For the latter we use, in each block, the lower Hashin-Shtrikman bound $\kappa_{LHS}$ (if $\kappa_M\leq\kappa_I$) or the upper Hashin-Shtrikman bound $\kappa_{UHS}$ (if $\kappa_M\geq\kappa_I$) \cite[p.557 and p.406]{torquato2013} computed with the sample-dependent volume fractions in that block. Denoting by $\phi_M$ and $\phi_I$ the (sample-dependent) volume fractions for matrix and inclusions in one block, $\kappa_{LHS}$ and $\kappa_{UHS}$ are given by
\begin{linenomath}
\begin{equation}\label{eq:HSbounds}
\kappa_{\rm HSL}  =  \kappa_{\rm arith} - \frac{(\kappa_M - \kappa_I)^2 \phi_M \phi_I}
{\kappa_I\phi_M + \kappa_M\phi_I + \kappa_{\rm min}},\quad 
\kappa_{\rm HSU}  =   \kappa_{\rm arith} - \frac{(\kappa_M - \kappa_I)^2 \phi_M \phi_I}
{\kappa_I\phi_M + \kappa_M\phi_I + \kappa_{\rm max}},
\end{equation} 
\end{linenomath}
with
$\kappa_{\rm arith} = \phi_M\kappa_M + \phi_I\kappa_I$,
$\kappa_{\rm min} = {\rm min} \{ \kappa_M, \kappa_I \}$, and
$\kappa_{\rm max} = {\rm max} \{ \kappa_M, \kappa_I \}$. 
We use blockwise homogenized coefficients instead of globally homogenized coefficients because
 they provide better variance reduction for the MLMC algorithm. The coefficient $\kappa_{eff}(\omega,\cdot)$, $\omega\in\Omega$ has to be understood as a \textsl{surrogate} for the actual, deterministic, homogenized coefficient.}

\prev{The} homogenized problem \eqref{eq:surg_forward_adj} involves the form
\begin{linenomath}
\begin{equation}
B_0(u_0,v) = \bbE\left[  \int_D \kappa_{0} \nabla u_0 \cdot \nabla v\dee\bfx \right]
\end{equation}
\end{linenomath}
and the error estimator \eqref{eq:eta_est} is 
\begin{linenomath}
\begin{equation}\label{eq:eta_est_heat}
 \eta_{est} = \mathbb{E}\left[-\int_{D\setminus D_{fine}} \kappa_0\left(1-\frac{\kappa_0}{\kappa}\right)\nabla u_0\cdot \nabla v_0 \dee \bfx\right].
\end{equation}
\end{linenomath}

It remains to define the subdomain $D_{fine}$ where the microstructure is resolved. \rev{This can be constructed adaptively using local error indicators. Namely, for each block $B_k$, $k=1,\ldots,K_b$ the quantity
\begin{linenomath}
\begin{equation}\label{eq:Nk}
\eta_{est,k} := \mathbb{E}\left[-\int_{B_k} \kappa_0\left(1-\frac{\kappa_0}{\kappa}\right)\nabla u_0\cdot \nabla v_0 \dee\bfx\right]
\end{equation}
\end{linenomath}
is an indicator for the contribution of the block $B_k$ to the total modeling error $\eta_{est}$, \prev{and $\sum_{k=1}^{K_b} \eta_{est,k} = \eta_{est}$.}
}
\prev{Similarly to \cite{odenvemaganti2000,Romkesetal}, the local error indicators can be used to guide the adaptive selection of surrogate models, as described in Subsection \ref{sssect:modelsel} in the MLMC framework.}

For later application of MLMC, we not only need the model fulfilling the required tolerance, but also a sequence of coarser models. Therefore, to define the coarsest possible model, before the blockwise homogenized one with $D_{fine}=\varnothing$, we consider the model with $\kappa_0(\omega)\equiv \kappa_{eff}(\omega)$, $\omega\in\Omega$, with $\kappa_{eff}(\omega)$ constant over $D$ and coinciding with the lower (resp. upper) Hashin-Shtrikman bound computed with the volume fractions over the whole domain.

\begin{rmk}[Cost of globally homogenized model]\label{rmk:scalarg}
 \prev{Since the globally homogenized coefficient is constant over the whole domain and we consider linear functionals of the solution, the cost of computing one realization of the QoI on the corsest model is $\mathcal{O}(1)$. Indeed, it is sufficient to solve the forward problem once to obtain the value $q_{fix}$ for the QoI with a fixed coefficient $\kappa_{fix}$, and any other realization with coefficient $\kappa_{eff}(\omega)$ can be obtained by a scaling of $q_{fix}$.} In particular, for \eqref{eq:sur_qoi1}, $q(\omega)=q_{fix}\frac{\kappa_{fix}}{\kappa_{eff}(\omega)}$, $\omega\in\Omega$.
\end{rmk}

\begin{rmk}[Additional surrogate models at coarse scale]\label{rmk:surrheat} \prev{If all blocks in which we partition the domain are quadrilaterals, then we can consider them as a mesh on the coarse scale. Namely, before resolving the fine scale in some blocks (but after the model with globally homogenized coefficient), we can consider a nested sequence of quadrilateral meshes, the finest one corresponding to the block partition, and associate to each of these a model where the (sample-dependent) elementwise homogenized coefficient is used. This procedure is clarified in Subsection \ref{ssec:numexp_heat}.}
\end{rmk}

\begin{rmk}[Possible extensions] \prev{The construction of surrogate models presented in this section can be coupled to other model reduction techniques available in the literature. In particular, in the models on the coarse scale, where the piecewise homogenized coefficient is used, the affine decomposition of the bilinear form would allow for an efficient application of the reduced basis method \cite{RB}. When the microstructure is resolved in some blocks, then application of the reduced basis method would require an empirical interpolation \cite{RBempirical}, which is potentially very expensive and does not always ensure a reduced basis construction which is sufficiently accurate.}
\end{rmk}

\subsection{Extension to Elastostatics of Random Materials}\label{sec:ela}\label{ssec:ela}

While algebraically more tedious, the full theory and adaptive algorithms described earlier are readily extendable to linear elasticity problems.
Dealing then with the vector-valued displacement field $\bfu = \bfu(\omega,\bfx), \bfu:\Omega \times D \rightarrow \mathbb{R}^{\prev{d}}$, \prev{$d=1,2,3$} ($\bfu \in \prev{\mathbf{H}} := (L^2(\Omega,\mathbb{P})\times H^1_{\Gamma_D}(D))^{\prev{d}}$), the model forward problem for a two-phase random elastic material with isotropic phases is of the form, for a.e. $\omega \in \Omega$, 
\begin{linenomath}
\begin{equation}\label{eq:elastic_forward}
\left.
\arraycolsep=1.4pt\def\arraystretch{1.4}
\begin{array}{rcll}
\nabla\cdot \bfT(\omega,\mathbf{u}) & = & \bff (\omega,\mathbf{x}), & \bfx \in D, \\
\bfT(\omega,\mathbf{u}) \bfn & = & \boldsymbol{\sigma}(\omega,\mathbf{x}), & \bfx \in \Gamma_N, \\
\bfu(\omega,\mathbf{x}) & = & 0, & \bfx \in \Gamma_D,
\end{array}
\right\}
\end{equation}
\end{linenomath}
with $\bfT$ the Cauchy stress, given by
\begin{linenomath}
\begin{equation}
\left.
\begin{aligned}
\bfT(\omega,\mathbf{u}) & =  2 \mu(\omega,\mathbf{x})\boldsymbol{\varepsilon}(\bfu) +  \lambda (\omega,\mathbf{x})\text{tr}(\boldsymbol{\varepsilon}(\bfu))\bfI, \\
\boldsymbol{\varepsilon}(\mathbf{u}) & =  \frac{1}{2} \left( \nabla \bfu + \left(\nabla\bfu\right)^{\top}\right).
\end{aligned}
\right\}
\end{equation}
\end{linenomath}
Here $\lambda$ and $\mu$ are the Lam\'e constants for the materials, which take on different values on each phase, $M$ and $I$, i.e., in analogy with \eqref{eq:kappaIM}, 
\begin{linenomath}
\begin{equation}
(\lambda, \mu) (\omega,\mathbf{x}) = (\lambda, \mu)_M\mathcal{X}_M(\omega,\mathbf{x}) + (\lambda, \mu)_I\mathcal{X}_I(\omega,\mathbf{x}), \qquad \omega\in\Omega, \quad \mathbf{x}\in\bar{D}.
\end{equation}
\end{linenomath}
\prev{The elasticity problem has the same structure as \eqref{eq:strongstochastic}, with the gradient replaced by its symmetric part and $\bfA(\omega,\bfx)$ a fourth-order tensor with components $A_{ijkl}=\lambda(\omega,\bfx)\delta_{ij}\delta_{kl} + \mu\left(\delta_{ik}\delta_{jl}+\delta_{il}\delta_{jk}\right)$, $1\leq i,j,k,l\leq d$ (where, for any $i,j=1,\ldots,d$, $\delta_{ij}$ is the Kronecker delta).}
%

The weak or variational form of \eqref{eq:elastic_forward} can be written
\begin{linenomath}
\begin{equation}
\left.
\begin{aligned}
\mbox{Find } \bfu \in \bfH &  \mbox{ such that for all} \; \bfv \in \bfH\\
B(\bfu, \bfv) &:=  F(\bfv) \\
B(\bfu, \bfv) &:= \int_\Omega \int_D 2\mu(\omega,\mathbf{x}) \boldsymbol{\varepsilon}(\bfu) : \boldsymbol{\varepsilon}(\bfv) + \lambda(\omega,\mathbf{x})\text{ div}\,\bfu\text{ div}\,\bfv \dee\bfx \dee\mathbb{P}(\omega) \\
F(\bfv) &:=  \int_\Omega \int_D \bff(\omega,\mathbf{x}) \cdot \bfv(\omega,\mathbf{x})  \dee\bfx \dee\mathbb{P}(\omega) +  \int_\Omega \int_{\Gamma_N} \boldsymbol{\sigma}(\omega,\mathbf{x}) \cdot \bfv(\omega,\mathbf{x})   \dee\bfx \dee\mathbb{P}(\omega),
\end{aligned}
\right\}
\end{equation}
\end{linenomath}
with an analogous adjoint problem. Among QoIs are the ensemble average of strains in a subdomain of \prev{volume} $A$ or over an inclusion, e.g. $Q(\bfu) = \mathbb{E}\left[ |A|^{-1} \int_A \varepsilon_{11}(\bfu(\cdot,\mathbf{x})) \dee\bfx\right]$.

For the surrogate models, we use the same construction as in Section \ref{sec:heat}, where now the Hashin-Shtrikman bounds on the conductivity are replaced by the Hashin-Shtrikman bounds on the bulk modulus $K$ and shear modulus $\mu$ \cite[p.570]{torquato2013},
\begin{linenomath}
\prev{\begin{align}
  K_{LHS} &= K_{arith} - \frac{(K_M-K_I)^2\phi_I\phi_M}{K_I\phi_M+K_M\phi_I + \frac{2(d-1)}{d}\mu_{\min}}, &  \mu_{LHS} &= \mu_{arith} - \frac{(\mu_M-\mu_I)^2\phi_I\phi_M}{\mu_I\phi_M+\mu_M\phi_I + H_{\min}},\label{eq:hsela1}\\
  K_{UHS} &= K_{arith} - \frac{(K_M-K_I)^2\phi_I\phi_M}{K_I\phi_M+K_M\phi_I + \frac{2(d-1)}{d}\mu_{\max}}, &
  \mu_{UHS} &= \mu_{arith} - \frac{(\mu_M-\mu_I)^2\phi_I\phi_M}{\mu_I\phi_M+\mu_M\phi_I + H_{\max}},
\label{eq:hsela2}
\end{align}}
\end{linenomath}
\prev{with $K_{arith}=K_I\phi_I+K_M\phi_M$, $K_{\min}=\min\left\{K_I,K_M\right\}$, $K_{\max}=\max\left\{K_I,K_M\right\}$ (similarly for $\mu$), $d$ the spatial dimension and}
\begin{linenomath}
\prev{\begin{equation*}
H_{\min} = \mu_{\min}\frac{dK_{\min}/2 + (d+1)(d-2)\mu_{\min}/d}{K_{\min} + 2\mu_{\min}} \qquad\quad H_{\max} = \mu_{\max}\frac{dK_{\max}/2 + (d+1)(d-2)\mu_{\max}/d}{K_{\max} + 2\mu_{\max}}.
\end{equation*}}
\end{linenomath}
\prev{From \eqref{eq:hsela1} and \eqref{eq:hsela2}, the corresponding values for $\lambda$ can be derived by the relationship $\lambda=K-\frac{2}{3}\mu$. }

\begin{rmk}[Comment on Remark \ref{rmk:scalarg} for elasticity]\label{rmk:noscalarg} \prev{In case of homogeneous materials having the same Poisson ratio but different Young's moduli, a similar scaling argument as in Remark \ref{rmk:scalarg} holds. However, even if the two materials in the random heterogeneous medium have the same Poisson ratio, the Hashin-Shtrikman bounds lead to an effective Poisson ratio which is in general different from the one of the two materials. Due to the nonlinear formula \eqref{eq:hsela1} for $K_{LHS}$ and $\mu_{LHS}$, $\lambda_{LHS}$ and $\mu_{LHS}$ scale differently for different volume fractions, and thus we cannot employ such a simple scaling argument as in Remark \ref{rmk:scalarg} (same observation holds when formula \eqref{eq:hsela2} is used).}
\end{rmk}

\section{Multilevel Monte Carlo on a Sequence of Surrogate Models}\label{sec6:MLMC}

In this section, we describe a MLMC strategy for the efficient approximation of $Q=\mathbb{E}[q(u)]$, $Q: H\rightarrow\mathbb{R}$ being the localized QoI and $u$ the solution to \eqref{eq:weakstochastic}. \prev{After a brief review of MLMC in Subsection \ref{ssec:MLMCintro}, we illustrate, in Subsection \ref{ssect:MLMCstrategy}, how the key features of MLMC can be combined with the model adaptive strategy to improve the efficiency of Monte Carlo sampling for local QoIs.}

\subsection{The Multilevel Monte Carlo method}\label{ssec:MLMCintro}
\prev{Here we follow the presentations in \cite{collier2015contMLMC,giles_2015}.} Let $q$ be the functional of interest, which is, in general, not accessible from computations, and let $\left(q_l\right)_{l=1}^{L}$ be a sequence of approximations to $q$, ordered from the least to the most accurate. Furthermore, we denote $q_0:=0$. The MLMC method exploits the identity
\begin{linenomath}
\begin{equation}\label{eq:mlmcid}
\mathbb{E}[q_L] = \sum_{l=1}^L \mathbb{E}\left[Y_l\right],\quad \text{with }Y_l:=q_l-q_{l-1},\, l=1,\ldots,L,
\end{equation}
\end{linenomath}
to estimate $Q=\mathbb{E}[q]$ by
\begin{linenomath}
\begin{equation}\label{eq:mlmcest1}
E^L[q_L] := \sum_{l=1}^L E_{M_l}\left[Y_l\right],
\end{equation}
\end{linenomath}
where
\begin{linenomath}
\begin{equation}\label{eq:mlmcest2}
E_{M_l}\left[Y_l\right]:=\frac{1}{M_l}\sum_{i=1}^{M_l} Y_l^i, \text{ for }l=1,\ldots,L.
\end{equation}
\end{linenomath}
In the equation above, for each $l=1,\ldots,L$, $\left\{Y_l^i\right\}_{i=1}^{M_l}$ denote $M_l\in\mathbb{N}$ i.i.d. (independent, identically distributed) samples of $Y_l$. 

For every realization of \eqref{eq:mlmcest1}, the error can be decomposed as
\begin{linenomath}
\begin{equation}\label{eq:mlmcerr}
\left|\mathbb{E}[q]-E^L[q_L]\right|\leq \left|\mathbb{E}[q-q_L]\right| + \left|\mathbb{E}[q_L]-E^L[q_L]\right|.
\end{equation}
\end{linenomath}
The first term on the right-hand side is deterministic, and it is called \textsl{bias error}. The second summand is a random variable, \prev{corresponding} to the \textsl{statistical error}. To compute the MLMC estimator up to a tolerance $TOL$, we can decompose \prev{the latter} into a bias tolerance $TOL_{bias}$ and a statistical error tolerance $TOL_{stat}$ \cite{collier2015contMLMC,hoel2014adaptiveMLMC}:
\begin{linenomath}
\begin{equation}\label{eq:tolsplit}
TOL = \underbrace{(1-\vartheta)TOL}_{:=TOL_{bias}} + \underbrace{\vartheta TOL}_{:=TOL_{stat}},
\end{equation}
\end{linenomath}
for some $\vartheta\in(0,1]$. \prev{For a fixed $\vartheta$, the bias tolerance is imposed by selecting the approximation $q_L$ such that}
\begin{linenomath}
\begin{equation}\label{eq:bias}
\left|\mathbb{E}[q-q_L]\right|\leq TOL_{bias}.
\end{equation}
\end{linenomath}
As described in the next subsection, this means choosing a sufficiently good surrogate model.

Due to its stochastic nature, there are multiple possibilities of prescribing a tolerance on \prev{the stochastic error}. \prev{A popular choice (see e.g. \cite{giles2008,giles_2015,cliffe2011}) is to do it} in the mean squared sense, requiring that 
\begin{linenomath}
\begin{equation}\label{eq:stat}
\mathbb{E}\left[\left|\mathbb{E}[q_L]-E^L[q_L]\right|^2\right]\leq TOL_{stat}^2.
\end{equation}
\end{linenomath}
This, together with \eqref{eq:bias}, ensures that
\begin{linenomath}
\begin{equation}\label{eq:mse}
\lVert\mathbb{E}[q]-E^L[q_L]\rVert_{L^2(\Omega,\mathbb{R})}\leq TOL.
\end{equation}
\end{linenomath}
The requirement in \eqref{eq:stat} is adopted in this paper. Alternatively, \prev{the statistical tolerance can be imposed in a probabilistic sense}, cf. \cite{collier2015contMLMC,hoel2014adaptiveMLMC}.

 The statistical error is controlled by 
 \begin{linenomath}
\begin{equation}
 \mathbb{E}\left[\left|\mathbb{E}[q_L]-E^L[q_L]\right|^2\right]\leq \sum_{l=1}^L \frac{\mathbb{V}[q_l-q_{l-1}]}{M_l}
\end{equation}
\end{linenomath}
(where $\mathbb{V}[\cdot]$ denotes the variance) \cite{giles_2015}. The estimated total work of the MLMC estimator is then \cite{giles_2015}
\begin{linenomath}
\begin{equation}\label{eq:cost}
 \mathbb{E}[W] = TOL_{stat}^{-2}\sum_{l=1}^L \sqrt{V_l W_l},
\end{equation}
\end{linenomath}
where $V_l=\mathbb{V}[Y_l]$, and $W_l$ is the estimated cost of one realization of $Y_l$, for $l=1,\ldots,L$. In this respect, the MLMC method can be seen as a variance reduction technique \cite{giles_2015, botevvariance}: the sequence $\left(q_l\right)_{l=1}^L$ must be such that, as $l$ increases and the cost of generating a sample increases, the variance of $Y_l$ decreases. This makes it possible to evaluate many samples on coarser levels and fewer on the finer ones, \prev{reducing the computational cost compared to plain Monte Carlo}. However, the complexity theorem for MLMC \cite[Thm. 2.1]{giles_2015} clearly establishes that \prev{such} gain in efficiency \prev{relies} on a very delicate trade-off between computational cost and accuracy across the levels. \prev{When dealing with a sequence of surrogate models, we cannot rely on convergence rates to determine the distribution of work across the levels. Defining a geometric sequence of tolerances as in \cite{hoel2014adaptiveMLMC} is also not applicable, as a discrete set of models does not ensure that for each level of tolerance a surrogate model exists and, also if it does, the relationship between computational costs and tolerances among the levels may be highly not optimal. How to tackle this issue is described in the next subsections.}


\subsection{Model-based Multilevel Monte Carlo}\label{ssect:MLMCstrategy}

Our aim is to use the same estimator as in \eqref{eq:mlmcest1} with $q_l:=q(u_0^l)$, where $u_0^l$ is the solution to a surrogate model as from Section \ref{sec:surrpb}, and thus
\begin{equation}
Y_l=q(u_0^l)-q(u_0^{l-1}),
\end{equation}
 for $l=1,\ldots,L$ and $L\in\mathbb{N}$.

The strategy that we use can be subdivided into three steps:
\begin{enumerate}
 \item model selection,
 \item level selection,
 \item computation of MLMC estimator.
\end{enumerate}

The model selection procedure leverages the a posteriori error estimator to construct a sequence of surrogate models. Out of this sequence, the level selection step selects $L$ models as levels for the MLMC estimator. Finally, the MLMC estimate of $\mathbb{E}[q]$ as in \eqref{eq:mlmcest1}-\eqref{eq:mlmcest2} is computed.

\subsubsection{Model selection}\label{sssect:modelsel}

For a given bias tolerance $TOL_{bias}$, the aim of the model selection procedure is twofold: to select the cheapest possible model $\mathcal{M}_J$ such that 
\begin{equation}
 \mathbb{E}\left[\left|q(u_0^{(J)})-q(u)\right|\right]\leq TOL_{bias},
\end{equation}
and to construct a sequence of surrogate models  with tolerance larger than the bias tolerance as candidates for the MLMC levels. 

\prev{Ideally, also in the deterministic setting, one would start with the coarsest model and enrich it adding microscale features until the estimated error is lower than $TOL_{bias}$. However, the quantity $\eta_{est}$ in \eqref{eq:eta_est} is an error \textsl{indicator} and not a reliable error estimator (the same holds for other error indicators introduced e.g. in \cite{odenvemaganti2000,Romkesetal}). Therefore, we use the local error estimators in \eqref{eq:Nk} to guide the adaptive process, similarly to \cite{eigel2016adaptive, hoel2014adaptiveMLMC}, and a sample average of the exact error $Q(e_0^{j})= \mathbb{E}\left[\left|q(u_0^{(j)})-q(u)\right|\right]$, $j=0,\ldots,J$, to quantify the modeling error of a surrogate.}

\prev{Let us denote by $\mathcal{S}$ the sequence of surrogate models that we want to construct, including $\mathcal{M}_J$ as the last model. We denote by $\mathcal{M}_0$ the model with blockwise homogenized coefficient and we start with $\mathcal{S}=\left\{\mathcal{M}_0\right\}$. We proceed then with the following steps, presented in Algorithm \ref{alg:modelselect}:
\begin{enumerate}
\item set $j=0$ and draw $\hat{M}$ samples of the microstructure (this does not require necessarily to save $\hat{M}$ meshes, but only the parameters or the random fields fully describing it);
\item compute and save the realizations of the QoI on the fine scale for the $\hat{M}$ samples, $q^i(u(\cdot))$, $i=1,\ldots,\hat{M}$; we will need these to estimate the exact error and, as already mentioned, ideally this step can be avoided if one could use a reliable error estimator; moreover, in view of the level selection, compute and save the sample average $E_{\hat{M}}[W^{fine}]$, where $W^{fine}$ is the cost of computing one sample of $q(u)$;
\item compute and save the $\hat{M}$ realizations of the QoI for the surrogate, $q^i(u_0^j)$,  $i=1,\ldots,\hat{M}$, and for these compute the errors $|q^i(u)-q^i(u_0^j)|$, $i=1,\ldots,\hat{M}$; also, compute and save the sample average $E_{\hat{M}}[W^j]$, where $W^j$ is the cost of computing one sample of $q(u^j)$;
\item if the sample average $E_{\hat{M}}\left[|q(u)-q(u_0^j)|\right]$ is lower than the bias tolerance, then stop, otherwise $j=j+1$ and select as next surrogate the one where microscale features are added in the blocks $B_k$, $k=1,\ldots,K_b$, such that
\begin{linenomath}
\begin{equation}\label{eq:blocksel}
E_{\hat{M}}\left[|\eta_{est,k}^{(j)}|\right]\geq \gamma E_{\hat{M}}\left[|\eta_{est,\tilde{k}}^{(j)}|\right],
\end{equation}
\end{linenomath}
where $\tilde{k}$ is the block with maximum local error indicator (in absolute value) and $\eta_{est,k}^{(j)}$ is the local error indicator for the $k$-th block and the $j$-th surrogate model; $\gamma\in(0,1)$ is an a priori chosen parameter, and in our experiments we used $\gamma=0.5$;
\item repeat steps 3-4 until $TOL_{bias}$ is reached, and correct the initial bias tolerance, setting it to the estimated error of the last model selected;
\item return the QoI samples and average cost for each surrogate and for the fine scale model, together with the set of models and the bias tolerance.
\end{enumerate}
If the scaling argument of Remark \ref{rmk:scalarg} applies (e.g. for the heat equation and linear QoI), then, denoting by $\mathcal{M}_{-1}$ the model in which the globally homogenized coefficient is used, we start instead with $\mathcal{S}=\left\{\mathcal{M}_{-1}\right\}$, perform step 3 for this model and proceed to $\mathcal{M}_0$ if the bias tolerance is exceeded. Moreover, if also Remark \ref{rmk:surrheat} is applicable, then we start with $\mathcal{S}=\left\{\mathcal{M}_{-1}\right\}$, perform step 3, and, if the tolerance is exceeded, we select the model associated to the coarsest quadrilateral mesh; then step 3 is iterated, refining each time the mesh if the tolerance is not reached; the finest subdivision in quadrilaterals corresponds to the model $\mathcal{M}_0$, and from there we can proceed as in Algorithm \ref{alg:modelselect}.

We note that the same $\hat{M}$ samples are considered for all sample averages. This and the fact that we save the QoI realizations for each surrogate will make the level selection procedure very cheap \cite{vidal2016}. The parameter $\hat{M}$ must be taken sufficiently large to ensure reliable sample estimates. In principle, one could control the accuracy of the sample estimates as in \cite[Algorithm 1]{hoel2014adaptiveMLMC}.

It remains to determine how to choose the bias tolerance, that is $\vartheta$ in \eqref{eq:tolsplit}. We start with the common splitting $TOL_{bias}=\tfrac{1}{\sqrt{2}}TOL$. When the statistical tolerance is corrected in step 5 (line \ref{algline:correcttol} of Algorithm \ref{alg:modelselect}), then an unbalanced splitting is used that tries to reduce the computational cost of the multilevel estimator for the prescribed total tolerance.}

\prev{\begin{algorithm}[H]
\caption{Error Estimator-driven Model Selection}\label{alg:modelselect}
\begin{algorithmic}[1]
\Statex \textsl{Input:} $TOL_{bias}$, $\hat{M}$, $\gamma$.
\Statex \textsl{Output:} $\mathcal{S}$, $TOL_{bias}$, QoI samples and average cost for each model.
\State $\mathcal{S}=\left\{\mathcal{M}_0\right\}$, \rev{$D_{fine}=\varnothing$}, $j=0$, $E_{\hat{M}}[\eta^{(0)}]=\text{Inf}$
\State Compute and save $\left\{q^i(u)\right\}_{i=1}^{\hat{M}}$ and $E_{\hat{M}}[W^{fine}]$
\State Compute and save $\left\{q^i(u_0^0)\right\}_{i=1}^{\hat{M}}$ and $E_{\hat{M}}[W^0]$
\If {$E_{\hat{M}}\left[|q(u)-q(u_0^0)|\right]<TOL_{bias}$}
\State $TOL_{bias}=E_{\hat{M}}\left[|q(u)-q(u_0^0)|\right]$
\State \Return $\mathcal{S}$, $TOL_{bias}$, $\left\{q^i(u_0^0)\right\}_{i=1}^{\hat{M}}$, $E_{\hat{M}}[W^0]$, $\left\{q^i(u)\right\}_{i=1}^{\hat{M}}$, $E_{\hat{M}}[W^{fine}]$ 
\Else 
\State $j=1$
\EndIf
\While{$E_{\hat{M}}[|q(u)-q(u_0^{j-1})|]>TOL_{bias}$}
\State Compute and save $\left\{q^i(u_0^j)\right\}_{i=1}^{\hat{M}}$ and $E_{\hat{M}}[W^j]$
\State Compute $E_{\hat{M}}[|q(u)-q(u_0^j)|]$,  $E_{\hat{M}}[\eta^{(j)}_k]$ for $k=1,\ldots,K_b$\label{algline:sample}
\If {$E_{\hat{M}}[|q(u)-q(u_0^j)|] > TOL_{bias}$}\label{algline:checktol}
\State $D_{fine}=D_{fine}\cup \mathcal{B}^j$, with $\mathcal{B}^j$ set of blocks that satisfy \eqref{eq:blocksel} $\rightarrow$ model $\mathcal{M}_j$\label{algline:addblk}
\State $\mathcal{S} = \mathcal{S}\cup \left\{\mathcal{M}_j\right\}$, $j=j+1$
\algstore{alg1}
\end{algorithmic}
\end{algorithm}
\begin{algorithm}
\begin{algorithmic}
\algrestore{alg1}
\If {all blocks have been added}
\State $TOL_{bias}=0$, $J=j-1$, \textbf{break}
\EndIf
\Else
\State $J=j$, $TOL_{bias}=E_{\hat{M}}\left[|q(u)-q(u_0^J)|\right]$, \textbf{break}\label{algline:correcttol}
\EndIf
\EndWhile
\State \Return $\mathcal{S}$, $TOL_{bias}$, $\left\{q^i(u_0^j)\right\}_{i=1}^{\hat{M}}$ and $E_{\hat{M}}[W^j]$ for $j=0,\ldots,J$, $\left\{q^i(u)\right\}_{i=1}^{\hat{M}}$,  $E_{\hat{M}}[W^{fine}]$
\end{algorithmic}
\end{algorithm}}

\subsubsection{Level selection}

\prev{Given a number of levels $L$, the goal of the level selection step is to extract the subsequence $\mathcal{S}_L:=\left\{\mathcal{M}_{j_1},\ldots\mathcal{M}_{j_L}\right\}\subseteq \mathcal{S}$ that minimizes the cost of the MLMC estimator for the prescribed tolerance $TOL$. Here we assume $L$ to be chosen a priori.

To select the subsequence $\mathcal{S}_L$, we proceed as follows:
\begin{enumerate}
\item if the adaptive procedure in the model selection has not reached the fine scale model, then consider two alternatives: using $TOL_{bias}$ as from the output of Algorithm \ref{alg:modelselect} and $\mathcal{S}=\left\{\mathcal{M}_1,\mathcal{M}_0,\ldots,\mathcal{M}_J\right\}$, or using $TOL_{bias}=0$ and $\mathcal{S}=\left\{\mathcal{M}_1,\mathcal{M}_0,\ldots,\mathcal{M}_J, \mathcal{M}_{fine}\right\}$, where $\mathcal{M}_{fine}$ is the fine scale model; if instead the adaptive procedure reached the fine scale model, consider the only choice $\mathcal{S}=\left\{\mathcal{M}_1,\mathcal{M}_0,\ldots,\mathcal{M}_J\right\}$ and $TOL_{bias}=0$;
\item for each option for $\mathcal{S}$ (two options if the fine scale model was not reached, one otherwise), set the last element of $\mathcal{S}_L$ as the last element of $\mathcal{S}$ and, for selecting the other levels, perform an exhaustive search over all ordered $(L-1)$-tuples in $\mathcal{S}$: for each of these $(L-1)$-tuples, compute the estimated cost as from \eqref{eq:cost} and choose the $(L-1)$-tuple that minimizes it; to compute the sample variances and sample averages in \eqref{eq:cost}, use the QoI realizations and average cost for each model saved in the model selection procedure;
\item if the fine scale model was not reached in the model selection, compare the estimated MLMC costs for the $L$-tuple for $TOL_{bias}>0$ and for the $L$-tuple for $TOL_{bias}=0$, and select the case with lower cost.
\end{enumerate}}

\subsubsection{Computation of the Multilevel Monte Carlo estimator}

Using the sequence $\mathcal{S}_L$ from the model selection, the standard MLMC estimator \prev{as in \cite{giles2008,cliffe2011}} can be computed. \prev{To compute the distribution of samples per level, we can use the sample averages of the variances and costs from the level selection procedure.} The work \cite{collier2015contMLMC} provides a strategy to avoid sample averages, but the algorithm proposed in that paper relies heavily on availability of convergence rates and therefore cannot be used directly for model hierarchies. 


\section{Numerical Experiments}\label{sec:numexp}

In this section, we present the results for the algorithms and theory given thusfar applied to representative model problems in heat conduction and in elastostatics in random, two-phase media.

For all cases, we have tested the MLMC algorithm for $L=2$ and $L=3$ levels. In each test case, the model selection has been performed once for both $L=2$ and $L=3$. In Algorithm \ref{alg:modelselect}, we have used \prev{$\gamma=0.5$, while $\hat{M}$ will be specified for each case. The local error indicators have been computed as in \eqref{eq:Nk}.} 

For the spatial discretization, conforming linear finite elements have been used. In our test cases, we consider a length scale in the microstructure which is \prev{not smaller than} ten percent of the length scale in the homogenized part and therefore using non-matching grids in the two regions would have not brought significant cost savings. However, when the ratio between the scales is higher, then using a fixed non-matching grid in the homogenized region brings significant cost savings. The works $W_l$, $l=1,\ldots,L$, are computed as $N_{dofs,l} + N_{dofs,l-1}$, where $N_{dofs,l}$ is the number of degrees of freedom at level $l$. These correspond to the standardised costs of solving the linear algebraic system from the PDE discretization, assuming an algebraic multigrid solver. Here we assume that the cost of solving the system dominates on the cost of meshing one sample. \prev{In the convergence plots, showing the root mean square error (left-hand side of \eqref{eq:mse}) versus work, the work for MLMC includes the preprocessing costs of the model selection procedure (the cost of the level selection is negligible); for the latter, the cost of each sample is computed as twice the number of degrees of freedom in the surrogate to take into account the solution of the adjoint problem (for the fine scale model the cost is only once the number of degrees of freedom, as the adjoint is not needed). The implementation is based on the FEniCS library \cite{fenicsI,fenicsII}.}

\subsection{Heat Conduction in Random Media}\label{ssec:numexp_heat}

\begin{figure}
\centering
 \includegraphics[scale=0.2]{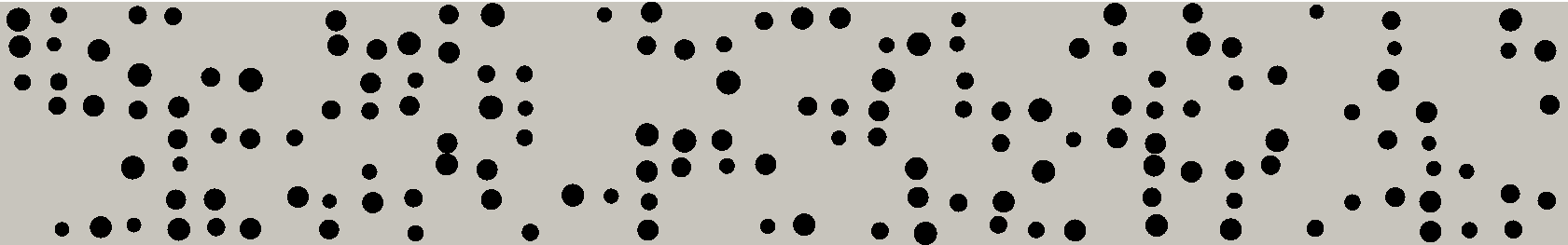}
\caption{\prev{One realization of} the geometry used for the heat conduction problem.}\label{fig:block}
\end{figure}

We consider a \prev{rectangular} domain $D=(0,2)\times (0,0.4)$ (m) with a random number of inclusions of circular shape, with stochastically perturbed radii and positions. \prev{A sample of the geometry is shown in Figure \ref{fig:block}. Remark \ref{rmk:scalarg} applies to the heat equation and linear QoIs, which we consider here. Having $D$ a tensor product structure, Remark \ref{rmk:surrheat} also applies, so that, before resolving the microstructure in some blocks, we consider the models on the coarse scale associated with the following diffusion coefficients and meshes:
\begin{itemize}
\item globally homogenized diffusion coefficient;
\item mesh of $5\times 1$  square elements (in $x-$ and $y-$ direction, respectively) with edge length $0.4$ and elementwise homogenized coefficient;
\item three uniform refinements of the $5\times 1$ mesh, with edge lengths $0.2$, $0.1$ and $0.05$, respectively; for each mesh, an elementwise homogenized coefficient is used.
\end{itemize}
The finest quadrilateral mesh corresponds to the subdivision of the domain into blocks as described in Subsection \ref{ssec:heat}.}

\prev{The generation of the inclusions is performed as follows:
\begin{itemize}
\item the number of inclusions per block is Bernoulli distributed $\sim \mathcal{B}er\left(\tfrac{1}{2}\right)$ (that is, each block has equal probability of containing or not one inclusion);
\item the $x$ and $y$ coordinates of the center of each inclusion are given by $x=x_c+\xi_x$ and $y=y_c+\xi_y$, with $\xi_x, \xi_y\sim\mathcal{U}\left(\left[-\frac{h}{8},\frac{h}{8}\right]\right)$(m), $h=0.05$(m);
\item the nominal radius of an inclusion, set to $r=\frac{h}{4}$(m), $h=0.05$(m), is perturbed by a uniform random variable $\mathcal{U}\left(\left[-\frac{h}{16},\frac{h}{16}\right]\right)$(m);
\item all random variables involved are independent.
\end{itemize}}

\medskip
\prev{The bilinear and linear forms in the fine scale model are given by  \eqref{eq:bilinheat}-\eqref{eq:linheat}, where $D = (0,2)\times (0,0.4)$,} the boundary $\Gamma_N$ coincides with the top and bottom boundaries of $D$, and $\Gamma_D$ are the left and right boundaries. The diffusion coefficient $\kappa$ is given by
\begin{linenomath}
\begin{equation*}
\kappa(\omega,x) = \begin{cases}
10,000 \; (W(mK)^{-1}) & \text{in the inclusions},\\
100 \; (W(mK)^{-1}) & \text{in the matrix},
\end{cases}
\end{equation*}
\end{linenomath}
and the applied flux $\sigma$ by
\begin{linenomath}
\prev{\begin{equation*}
\sigma(x) = \begin{cases}
1,600 \; (W⋅m^{-2})& \text{on the top boundary},\\
0 & \text{on the bottom boundary}.
\end{cases}
\end{equation*}}
\end{linenomath}

Since, as already mentioned, we focus on the modeling error, we need to use sufficiently fine meshes to ensure that the discretization error is negligible with respect to the modeling error. In this respect, we highlight that the performance of the model-based MLMC algorithm depends on the ratio between the mesh size in the homogenized region and the mesh size in the blocks which are refined, rather than the values of the mesh sizes themselves. \prev{Apart from the quadrilateral meshes on the coarse scale, for the other models (where the microscale is resolved in some blocks)} we use unstructured, conformal grids with a mesh size of \prev{$0.004$} in the refined blocks and \prev{$0.038$} in the homogenized region, leading to about \prev{$74050$} degrees of freedom for the \rev{fine scale} model (the actual number depends on the sample). \prev{As QoI we consider the average $y$-component of the gradient of the solution, denoted by $(\nabla u)_y$, in the block $B_q:=[0.1,0.2]\times [0.2,0.3]$:}
\begin{linenomath}
\prev{\begin{equation*}
q_1(u) = \frac{1}{A_q}\int_{B_q} (\nabla u)_y(\mathbf{x})\,\textit{d}\mathbf{x},
\end{equation*}}
\end{linenomath}
\prev{$A_q=0.001$($\text{m}^2$).} The solution $q_{fix}$ of the globally homogenized model is computed with $\kappa_{fix}=100$ on \prev{the coarsest quadrilateral mesh}. 

\prev{We consider the tolerances $TOL_0=0.02, TOL_1=0.01, TOL_2=0.005, TOL_3=0.0025$ for $L=2$ and $L=3$. Algorithm \ref{alg:modelselect} has been run with $\hat{M}=180$.}

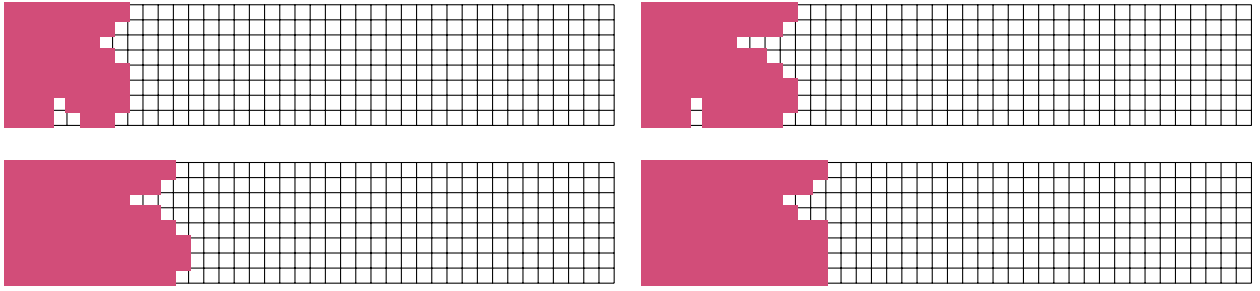
\begin{figure}
\centering
\begin{tikzpicture}[scale=0.2,every node/.style={minimum size=0.2cm-\pgflinewidth, outer sep=0pt}]
    \draw[step=1cm,color=black] (0,0) grid (40,8);
    \node[fill=purple!70] at (0.5,0.5) {};
    \node[fill=purple!70] at (1.5,0.5) {};
    \node[fill=purple!70] at (2.5,0.5) {};
    \node[fill=purple!70] at (5.5,0.5) {};
    \node[fill=purple!70] at (6.5,0.5) {};
    \node[fill=purple!70] at (0.5,1.5) {};
    \node[fill=purple!70] at (1.5,1.5) {};
    \node[fill=purple!70] at (2.5,1.5) {};
    \node[fill=purple!70] at (4.5,1.5) {};
    \node[fill=purple!70] at (5.5,1.5) {};
    \node[fill=purple!70] at (6.5,1.5) {};
    \node[fill=purple!70] at (7.5,1.5) {};
    \node[fill=purple!70] at (0.5,2.5) {};
    \node[fill=purple!70] at (1.5,2.5) {};
    \node[fill=purple!70] at (2.5,2.5) {};
    \node[fill=purple!70] at (3.5,2.5) {};
    \node[fill=purple!70] at (4.5,2.5) {};
    \node[fill=purple!70] at (5.5,2.5) {};
    \node[fill=purple!70] at (6.5,2.5) {};
    \node[fill=purple!70] at (7.5,2.5) {};
    \node[fill=purple!70] at (0.5,3.5) {};
    \node[fill=purple!70] at (1.5,3.5) {};
    \node[fill=purple!70] at (2.5,3.5) {};
    \node[fill=purple!70] at (3.5,3.5) {};
    \node[fill=purple!70] at (4.5,3.5) {};
    \node[fill=purple!70] at (5.5,3.5) {};
    \node[fill=purple!70] at (6.5,3.5) {};
    \node[fill=purple!70] at (7.5,3.5) {};
    \node[fill=purple!70] at (0.5,4.5) {};
    \node[fill=purple!70] at (1.5,4.5) {};
    \node[fill=purple!70] at (2.5,4.5) {};
    \node[fill=purple!70] at (3.5,4.5) {};
    \node[fill=purple!70] at (4.5,4.5) {};
    \node[fill=purple!70] at (5.5,4.5) {};
    \node[fill=purple!70] at (6.5,4.5) {};
    \node[fill=purple!70] at (0.5,5.5) {};
    \node[fill=purple!70] at (1.5,5.5) {};
    \node[fill=purple!70] at (2.5,5.5) {};
    \node[fill=purple!70] at (3.5,5.5) {};
    \node[fill=purple!70] at (4.5,5.5) {};
    \node[fill=purple!70] at (5.5,5.5) {};
    \node[fill=purple!70] at (0.5,6.5) {};
    \node[fill=purple!70] at (1.5,6.5) {};
    \node[fill=purple!70] at (2.5,6.5) {};
    \node[fill=purple!70] at (3.5,6.5) {};
    \node[fill=purple!70] at (4.5,6.5) {};
    \node[fill=purple!70] at (5.5,6.5) {};
    \node[fill=purple!70] at (6.5,6.5) {};
    \node[fill=purple!70] at (0.5,7.5) {};
    \node[fill=purple!70] at (1.5,7.5) {};
    \node[fill=purple!70] at (2.5,7.5) {};
    \node[fill=purple!70] at (3.5,7.5) {};
    \node[fill=purple!70] at (4.5,7.5) {};
    \node[fill=purple!70] at (5.5,7.5) {};
    \node[fill=purple!70] at (6.5,7.5) {};
    \node[fill=purple!70] at (7.5,7.5) {};
\end{tikzpicture}\hspace{0.35cm}\begin{tikzpicture}[scale=0.2,every node/.style={minimum size=0.2cm-\pgflinewidth, outer sep=0pt}]
    \draw[step=1cm,color=black] (0,0) grid (40,8);
    \node[fill=purple!70] at (0.5,0.5) {};
    \node[fill=purple!70] at (1.5,0.5) {};
    \node[fill=purple!70] at (2.5,0.5) {};
    \node[fill=purple!70] at (4.5,0.5) {};
    \node[fill=purple!70] at (5.5,0.5) {};
    \node[fill=purple!70] at (6.5,0.5) {};
    \node[fill=purple!70] at (7.5,0.5) {};
    \node[fill=purple!70] at (8.5,0.5) {};
    \node[fill=purple!70] at (0.5,1.5) {};
    \node[fill=purple!70] at (1.5,1.5) {};
    \node[fill=purple!70] at (2.5,1.5) {};
    \node[fill=purple!70] at (4.5,1.5) {};
    \node[fill=purple!70] at (5.5,1.5) {};
    \node[fill=purple!70] at (6.5,1.5) {};
    \node[fill=purple!70] at (7.5,1.5) {};
    \node[fill=purple!70] at (8.5,1.5) {};
    \node[fill=purple!70] at (9.5,1.5) {};
    \node[fill=purple!70] at (0.5,2.5) {};
    \node[fill=purple!70] at (1.5,2.5) {};
    \node[fill=purple!70] at (2.5,2.5) {};
    \node[fill=purple!70] at (3.5,2.5) {};
    \node[fill=purple!70] at (4.5,2.5) {};
    \node[fill=purple!70] at (5.5,2.5) {};
    \node[fill=purple!70] at (6.5,2.5) {};
    \node[fill=purple!70] at (7.5,2.5) {};
    \node[fill=purple!70] at (8.5,2.5) {};
    \node[fill=purple!70] at (9.5,2.5) {};
    \node[fill=purple!70] at (0.5,3.5) {};
    \node[fill=purple!70] at (1.5,3.5) {};
    \node[fill=purple!70] at (2.5,3.5) {};
    \node[fill=purple!70] at (3.5,3.5) {};
    \node[fill=purple!70] at (4.5,3.5) {};
    \node[fill=purple!70] at (5.5,3.5) {};
    \node[fill=purple!70] at (6.5,3.5) {};
    \node[fill=purple!70] at (7.5,3.5) {};
    \node[fill=purple!70] at (8.5,3.5) {};
    \node[fill=purple!70] at (0.5,4.5) {};
    \node[fill=purple!70] at (1.5,4.5) {};
    \node[fill=purple!70] at (2.5,4.5) {};
    \node[fill=purple!70] at (3.5,4.5) {};
    \node[fill=purple!70] at (4.5,4.5) {};
    \node[fill=purple!70] at (5.5,4.5) {};
    \node[fill=purple!70] at (6.5,4.5) {};
    \node[fill=purple!70] at (7.5,4.5) {};
    \node[fill=purple!70] at (0.5,5.5) {};
    \node[fill=purple!70] at (1.5,5.5) {};
    \node[fill=purple!70] at (2.5,5.5) {};
    \node[fill=purple!70] at (3.5,5.5) {};
    \node[fill=purple!70] at (4.5,5.5) {};
    \node[fill=purple!70] at (5.5,5.5) {};
    \node[fill=purple!70] at (0.5,6.5) {};
    \node[fill=purple!70] at (1.5,6.5) {};
    \node[fill=purple!70] at (2.5,6.5) {};
    \node[fill=purple!70] at (3.5,6.5) {};
    \node[fill=purple!70] at (4.5,6.5) {};
    \node[fill=purple!70] at (5.5,6.5) {};
    \node[fill=purple!70] at (6.5,6.5) {};
    \node[fill=purple!70] at (7.5,6.5) {};
    \node[fill=purple!70] at (8.5,6.5) {};
    \node[fill=purple!70] at (0.5,7.5) {};
    \node[fill=purple!70] at (1.5,7.5) {};
    \node[fill=purple!70] at (2.5,7.5) {};
    \node[fill=purple!70] at (3.5,7.5) {};
    \node[fill=purple!70] at (4.5,7.5) {};
    \node[fill=purple!70] at (5.5,7.5) {};
    \node[fill=purple!70] at (6.5,7.5) {};
    \node[fill=purple!70] at (7.5,7.5) {};
    \node[fill=purple!70] at (8.5,7.5) {};
    \node[fill=purple!70] at (9.5,7.5) {};
\end{tikzpicture}

\vspace{0.4cm}
\begin{tikzpicture}[scale=0.2,every node/.style={minimum size=0.2cm-\pgflinewidth, outer sep=0pt}]
    \draw[step=1cm,color=black] (0,0) grid (40,8);
    \node[fill=purple!70] at (0.5,0.5) {};
    \node[fill=purple!70] at (1.5,0.5) {};
    \node[fill=purple!70] at (2.5,0.5) {};
    \node[fill=purple!70] at (3.5,0.5) {};
    \node[fill=purple!70] at (4.5,0.5) {};
    \node[fill=purple!70] at (5.5,0.5) {};
    \node[fill=purple!70] at (6.5,0.5) {};
    \node[fill=purple!70] at (7.5,0.5) {};
    \node[fill=purple!70] at (8.5,0.5) {};
    \node[fill=purple!70] at (9.5,0.5) {};
    \node[fill=purple!70] at (10.5,0.5) {};
    \node[fill=purple!70] at (0.5,1.5) {};
    \node[fill=purple!70] at (1.5,1.5) {};
    \node[fill=purple!70] at (2.5,1.5) {};
    \node[fill=purple!70] at (3.5,1.5) {};
    \node[fill=purple!70] at (4.5,1.5) {};
    \node[fill=purple!70] at (5.5,1.5) {};
    \node[fill=purple!70] at (6.5,1.5) {};
    \node[fill=purple!70] at (7.5,1.5) {};
    \node[fill=purple!70] at (8.5,1.5) {};
    \node[fill=purple!70] at (9.5,1.5) {};
    \node[fill=purple!70] at (10.5,1.5) {};
    \node[fill=purple!70] at (11.5,1.5) {};
    \node[fill=purple!70] at (0.5,2.5) {};
    \node[fill=purple!70] at (1.5,2.5) {};
    \node[fill=purple!70] at (2.5,2.5) {};
    \node[fill=purple!70] at (3.5,2.5) {};
    \node[fill=purple!70] at (4.5,2.5) {};
    \node[fill=purple!70] at (5.5,2.5) {};
    \node[fill=purple!70] at (6.5,2.5) {};
    \node[fill=purple!70] at (7.5,2.5) {};
    \node[fill=purple!70] at (8.5,2.5) {};
    \node[fill=purple!70] at (9.5,2.5) {};
    \node[fill=purple!70] at (10.5,2.5) {};
    \node[fill=purple!70] at (11.5,2.5) {};
    \node[fill=purple!70] at (0.5,3.5) {};
    \node[fill=purple!70] at (1.5,3.5) {};
    \node[fill=purple!70] at (2.5,3.5) {};
    \node[fill=purple!70] at (3.5,3.5) {};
    \node[fill=purple!70] at (4.5,3.5) {};
    \node[fill=purple!70] at (5.5,3.5) {};
    \node[fill=purple!70] at (6.5,3.5) {};
    \node[fill=purple!70] at (7.5,3.5) {};
    \node[fill=purple!70] at (8.5,3.5) {};
    \node[fill=purple!70] at (9.5,3.5) {};
    \node[fill=purple!70] at (10.5,3.5) {};
    \node[fill=purple!70] at (0.5,4.5) {};
    \node[fill=purple!70] at (1.5,4.5) {};
    \node[fill=purple!70] at (2.5,4.5) {};
    \node[fill=purple!70] at (3.5,4.5) {};
    \node[fill=purple!70] at (4.5,4.5) {};
    \node[fill=purple!70] at (5.5,4.5) {};
    \node[fill=purple!70] at (6.5,4.5) {};
    \node[fill=purple!70] at (7.5,4.5) {};
    \node[fill=purple!70] at (8.5,4.5) {};
    \node[fill=purple!70] at (9.5,4.5) {};
    \node[fill=purple!70] at (0.5,5.5) {};
    \node[fill=purple!70] at (1.5,5.5) {};
    \node[fill=purple!70] at (2.5,5.5) {};
    \node[fill=purple!70] at (3.5,5.5) {};
    \node[fill=purple!70] at (4.5,5.5) {};
    \node[fill=purple!70] at (5.5,5.5) {};
    \node[fill=purple!70] at (6.5,5.5) {};
    \node[fill=purple!70] at (7.5,5.5) {};
    \node[fill=purple!70] at (0.5,6.5) {};
    \node[fill=purple!70] at (1.5,6.5) {};
    \node[fill=purple!70] at (2.5,6.5) {};
    \node[fill=purple!70] at (3.5,6.5) {};
    \node[fill=purple!70] at (4.5,6.5) {};
    \node[fill=purple!70] at (5.5,6.5) {};
    \node[fill=purple!70] at (6.5,6.5) {};
    \node[fill=purple!70] at (7.5,6.5) {};
    \node[fill=purple!70] at (8.5,6.5) {};
    \node[fill=purple!70] at (9.5,6.5) {};
    \node[fill=purple!70] at (0.5,7.5) {};
    \node[fill=purple!70] at (1.5,7.5) {};
    \node[fill=purple!70] at (2.5,7.5) {};
    \node[fill=purple!70] at (3.5,7.5) {};
    \node[fill=purple!70] at (4.5,7.5) {};
    \node[fill=purple!70] at (5.5,7.5) {};
    \node[fill=purple!70] at (6.5,7.5) {};
    \node[fill=purple!70] at (7.5,7.5) {};
    \node[fill=purple!70] at (8.5,7.5) {};
    \node[fill=purple!70] at (9.5,7.5) {};
    \node[fill=purple!70] at (10.5,7.5) {};
\end{tikzpicture}\hspace{0.35cm}\begin{tikzpicture}[scale=0.2,every node/.style={minimum size=0.2cm-\pgflinewidth, outer sep=0pt}]
    \draw[step=1cm,color=black] (0,0) grid (40,8);
    \node[fill=purple!70] at (0.5,0.5) {};
    \node[fill=purple!70] at (1.5,0.5) {};
    \node[fill=purple!70] at (2.5,0.5) {};
    \node[fill=purple!70] at (3.5,0.5) {};
    \node[fill=purple!70] at (4.5,0.5) {};
    \node[fill=purple!70] at (5.5,0.5) {};
    \node[fill=purple!70] at (6.5,0.5) {};
    \node[fill=purple!70] at (7.5,0.5) {};
    \node[fill=purple!70] at (8.5,0.5) {};
    \node[fill=purple!70] at (9.5,0.5) {};
    \node[fill=purple!70] at (10.5,0.5) {};
    \node[fill=purple!70] at (11.5,0.5) {};
    \node[fill=purple!70] at (0.5,1.5) {};
    \node[fill=purple!70] at (1.5,1.5) {};
    \node[fill=purple!70] at (2.5,1.5) {};
    \node[fill=purple!70] at (3.5,1.5) {};
    \node[fill=purple!70] at (4.5,1.5) {};
    \node[fill=purple!70] at (5.5,1.5) {};
    \node[fill=purple!70] at (6.5,1.5) {};
    \node[fill=purple!70] at (7.5,1.5) {};
    \node[fill=purple!70] at (8.5,1.5) {};
    \node[fill=purple!70] at (9.5,1.5) {};
    \node[fill=purple!70] at (10.5,1.5) {};
    \node[fill=purple!70] at (11.5,1.5) {};
    \node[fill=purple!70] at (0.5,2.5) {};
    \node[fill=purple!70] at (1.5,2.5) {};
    \node[fill=purple!70] at (2.5,2.5) {};
    \node[fill=purple!70] at (3.5,2.5) {};
    \node[fill=purple!70] at (4.5,2.5) {};
    \node[fill=purple!70] at (5.5,2.5) {};
    \node[fill=purple!70] at (6.5,2.5) {};
    \node[fill=purple!70] at (7.5,2.5) {};
    \node[fill=purple!70] at (8.5,2.5) {};
    \node[fill=purple!70] at (9.5,2.5) {};
    \node[fill=purple!70] at (10.5,2.5) {};
    \node[fill=purple!70] at (11.5,2.5) {};
    \node[fill=purple!70] at (0.5,3.5) {};
    \node[fill=purple!70] at (1.5,3.5) {};
    \node[fill=purple!70] at (2.5,3.5) {};
    \node[fill=purple!70] at (3.5,3.5) {};
    \node[fill=purple!70] at (4.5,3.5) {};
    \node[fill=purple!70] at (5.5,3.5) {};
    \node[fill=purple!70] at (6.5,3.5) {};
    \node[fill=purple!70] at (7.5,3.5) {};
    \node[fill=purple!70] at (8.5,3.5) {};
    \node[fill=purple!70] at (9.5,3.5) {};
    \node[fill=purple!70] at (10.5,3.5) {};
    \node[fill=purple!70] at (11.5,3.5) {};
    \node[fill=purple!70] at (0.5,4.5) {};
    \node[fill=purple!70] at (1.5,4.5) {};
    \node[fill=purple!70] at (2.5,4.5) {};
    \node[fill=purple!70] at (3.5,4.5) {};
    \node[fill=purple!70] at (4.5,4.5) {};
    \node[fill=purple!70] at (5.5,4.5) {};
    \node[fill=purple!70] at (6.5,4.5) {};
    \node[fill=purple!70] at (7.5,4.5) {};
    \node[fill=purple!70] at (8.5,4.5) {};
    \node[fill=purple!70] at (9.5,4.5) {};
    \node[fill=purple!70] at (0.5,5.5) {};
    \node[fill=purple!70] at (1.5,5.5) {};
    \node[fill=purple!70] at (2.5,5.5) {};
    \node[fill=purple!70] at (3.5,5.5) {};
    \node[fill=purple!70] at (4.5,5.5) {};
    \node[fill=purple!70] at (5.5,5.5) {};
    \node[fill=purple!70] at (6.5,5.5) {};
    \node[fill=purple!70] at (7.5,5.5) {};
    \node[fill=purple!70] at (8.5,5.5) {};
    \node[fill=purple!70] at (0.5,6.5) {};
    \node[fill=purple!70] at (1.5,6.5) {};
    \node[fill=purple!70] at (2.5,6.5) {};
    \node[fill=purple!70] at (3.5,6.5) {};
    \node[fill=purple!70] at (4.5,6.5) {};
    \node[fill=purple!70] at (5.5,6.5) {};
    \node[fill=purple!70] at (6.5,6.5) {};
    \node[fill=purple!70] at (7.5,6.5) {};
    \node[fill=purple!70] at (8.5,6.5) {};
    \node[fill=purple!70] at (9.5,6.5) {};
    \node[fill=purple!70] at (10.5,6.5) {};
    \node[fill=purple!70] at (0.5,7.5) {};
    \node[fill=purple!70] at (1.5,7.5) {};
    \node[fill=purple!70] at (2.5,7.5) {};
    \node[fill=purple!70] at (3.5,7.5) {};
    \node[fill=purple!70] at (4.5,7.5) {};
    \node[fill=purple!70] at (5.5,7.5) {};
    \node[fill=purple!70] at (6.5,7.5) {};
    \node[fill=purple!70] at (7.5,7.5) {};
    \node[fill=purple!70] at (8.5,7.5) {};
    \node[fill=purple!70] at (9.5,7.5) {};
    \node[fill=purple!70] at (10.5,7.5) {};
    \node[fill=purple!70] at (11.5,7.5) {};
\end{tikzpicture}
\caption{QoI $Q_1$, $L=2$: models chosen as level $2$ for $TOL_0$ (top left), $TOL_1$ (top right), $TOL_2$ (bottom left) and $TOL_3$ (bottom right).}\label{fig:L2levelsQ1}
\end{figure}

\begin{figure}
\centering
\begin{tikzpicture}[scale=0.2,every node/.style={minimum size=0.2cm-\pgflinewidth, outer sep=0pt}]
    \draw[step=1cm,color=black] (0,0) grid (40,8);
\node[fill=purple!70] at (0.5,0.5) {};
    \node[fill=purple!70] at (1.5,0.5) {};
    \node[fill=purple!70] at (2.5,0.5) {};
    \node[fill=purple!70] at (5.5,0.5) {};
    \node[fill=purple!70] at (6.5,0.5) {};
    \node[fill=purple!70] at (0.5,1.5) {};
    \node[fill=purple!70] at (1.5,1.5) {};
    \node[fill=purple!70] at (2.5,1.5) {};
    \node[fill=purple!70] at (4.5,1.5) {};
    \node[fill=purple!70] at (5.5,1.5) {};
    \node[fill=purple!70] at (6.5,1.5) {};
    \node[fill=purple!70] at (7.5,1.5) {};
    \node[fill=purple!70] at (0.5,2.5) {};
    \node[fill=purple!70] at (1.5,2.5) {};
    \node[fill=purple!70] at (2.5,2.5) {};
    \node[fill=purple!70] at (3.5,2.5) {};
    \node[fill=purple!70] at (4.5,2.5) {};
    \node[fill=purple!70] at (5.5,2.5) {};
    \node[fill=purple!70] at (6.5,2.5) {};
    \node[fill=purple!70] at (7.5,2.5) {};
    \node[fill=purple!70] at (0.5,3.5) {};
    \node[fill=purple!70] at (1.5,3.5) {};
    \node[fill=purple!70] at (2.5,3.5) {};
    \node[fill=purple!70] at (3.5,3.5) {};
    \node[fill=purple!70] at (4.5,3.5) {};
    \node[fill=purple!70] at (5.5,3.5) {};
    \node[fill=purple!70] at (6.5,3.5) {};
    \node[fill=purple!70] at (7.5,3.5) {};
    \node[fill=purple!70] at (0.5,4.5) {};
    \node[fill=purple!70] at (1.5,4.5) {};
    \node[fill=purple!70] at (2.5,4.5) {};
    \node[fill=purple!70] at (3.5,4.5) {};
    \node[fill=purple!70] at (4.5,4.5) {};
    \node[fill=purple!70] at (5.5,4.5) {};
    \node[fill=purple!70] at (6.5,4.5) {};
    \node[fill=purple!70] at (0.5,5.5) {};
    \node[fill=purple!70] at (1.5,5.5) {};
    \node[fill=purple!70] at (2.5,5.5) {};
    \node[fill=purple!70] at (3.5,5.5) {};
    \node[fill=purple!70] at (4.5,5.5) {};
    \node[fill=purple!70] at (5.5,5.5) {};
    \node[fill=purple!70] at (0.5,6.5) {};
    \node[fill=purple!70] at (1.5,6.5) {};
    \node[fill=purple!70] at (2.5,6.5) {};
    \node[fill=purple!70] at (3.5,6.5) {};
    \node[fill=purple!70] at (4.5,6.5) {};
    \node[fill=purple!70] at (5.5,6.5) {};
    \node[fill=purple!70] at (6.5,6.5) {};
    \node[fill=purple!70] at (0.5,7.5) {};
    \node[fill=purple!70] at (1.5,7.5) {};
    \node[fill=purple!70] at (2.5,7.5) {};
    \node[fill=purple!70] at (3.5,7.5) {};
    \node[fill=purple!70] at (4.5,7.5) {};
    \node[fill=purple!70] at (5.5,7.5) {};
    \node[fill=purple!70] at (6.5,7.5) {};
\end{tikzpicture}\hspace{0.35cm}\begin{tikzpicture}[scale=0.2,every node/.style={minimum size=0.2cm-\pgflinewidth, outer sep=0pt}]
    \draw[step=1cm,color=black] (0,0) grid (40,8);
    \node[fill=purple!70] at (0.5,0.5) {};
    \node[fill=purple!70] at (1.5,0.5) {};
    \node[fill=purple!70] at (0.5,1.5) {};
    \node[fill=purple!70] at (1.5,1.5) {};
    \node[fill=purple!70] at (2.5,1.5) {};
    \node[fill=purple!70] at (0.5,2.5) {};
    \node[fill=purple!70] at (1.5,2.5) {};
    \node[fill=purple!70] at (2.5,2.5) {};
    \node[fill=purple!70] at (4.5,2.5) {};
    \node[fill=purple!70] at (5.5,2.5) {};
    \node[fill=purple!70] at (0.5,3.5) {};
    \node[fill=purple!70] at (1.5,3.5) {};
    \node[fill=purple!70] at (2.5,3.5) {};
    \node[fill=purple!70] at (3.5,3.5) {};
    \node[fill=purple!70] at (4.5,3.5) {};
    \node[fill=purple!70] at (5.5,3.5) {};
    \node[fill=purple!70] at (6.5,3.5) {};
    \node[fill=purple!70] at (0.5,4.5) {};
    \node[fill=purple!70] at (1.5,4.5) {};
    \node[fill=purple!70] at (2.5,4.5) {};
    \node[fill=purple!70] at (3.5,4.5) {};
    \node[fill=purple!70] at (4.5,4.5) {};
    \node[fill=purple!70] at (5.5,4.5) {};
    \node[fill=purple!70] at (0.5,5.5) {};
    \node[fill=purple!70] at (1.5,5.5) {};
    \node[fill=purple!70] at (2.5,5.5) {};
    \node[fill=purple!70] at (3.5,5.5) {};
    \node[fill=purple!70] at (4.5,5.5) {};
    \node[fill=purple!70] at (0.5,6.5) {};
    \node[fill=purple!70] at (1.5,6.5) {};
    \node[fill=purple!70] at (2.5,6.5) {};
    \node[fill=purple!70] at (3.5,6.5) {};
    \node[fill=purple!70] at (4.5,6.5) {};
    \node[fill=purple!70] at (5.5,6.5) {};
    \node[fill=purple!70] at (0.5,7.5) {};
    \node[fill=purple!70] at (1.5,7.5) {};
    \node[fill=purple!70] at (2.5,7.5) {};
    \node[fill=purple!70] at (3.5,7.5) {};
    \node[fill=purple!70] at (4.5,7.5) {};
    \node[fill=purple!70] at (5.5,7.5) {};
    \node[fill=purple!70] at (6.5,7.5) {};
\end{tikzpicture}

\vspace{0.4cm}
\begin{tikzpicture}[scale=0.2,every node/.style={minimum size=0.2cm-\pgflinewidth, outer sep=0pt}]
    \draw[step=1cm,color=black] (0,0) grid (40,8);
    \node[fill=purple!70] at (0.5,0.5) {};
    \node[fill=purple!70] at (1.5,0.5) {};
    \node[fill=purple!70] at (2.5,0.5) {};
    \node[fill=purple!70] at (4.5,0.5) {};
    \node[fill=purple!70] at (5.5,0.5) {};
    \node[fill=purple!70] at (6.5,0.5) {};
    \node[fill=purple!70] at (7.5,0.5) {};
    \node[fill=purple!70] at (0.5,1.5) {};
    \node[fill=purple!70] at (1.5,1.5) {};
    \node[fill=purple!70] at (2.5,1.5) {};
    \node[fill=purple!70] at (4.5,1.5) {};
    \node[fill=purple!70] at (5.5,1.5) {};
    \node[fill=purple!70] at (6.5,1.5) {};
    \node[fill=purple!70] at (7.5,1.5) {};
    \node[fill=purple!70] at (0.5,2.5) {};
    \node[fill=purple!70] at (1.5,2.5) {};
    \node[fill=purple!70] at (2.5,2.5) {};
    \node[fill=purple!70] at (3.5,2.5) {};
    \node[fill=purple!70] at (4.5,2.5) {};
    \node[fill=purple!70] at (5.5,2.5) {};
    \node[fill=purple!70] at (6.5,2.5) {};
    \node[fill=purple!70] at (7.5,2.5) {};
    \node[fill=purple!70] at (0.5,3.5) {};
    \node[fill=purple!70] at (1.5,3.5) {};
    \node[fill=purple!70] at (2.5,3.5) {};
    \node[fill=purple!70] at (3.5,3.5) {};
    \node[fill=purple!70] at (4.5,3.5) {};
    \node[fill=purple!70] at (5.5,3.5) {};
    \node[fill=purple!70] at (6.5,3.5) {};
    \node[fill=purple!70] at (7.5,3.5) {};
    \node[fill=purple!70] at (0.5,4.5) {};
    \node[fill=purple!70] at (1.5,4.5) {};
    \node[fill=purple!70] at (2.5,4.5) {};
    \node[fill=purple!70] at (3.5,4.5) {};
    \node[fill=purple!70] at (4.5,4.5) {};
    \node[fill=purple!70] at (5.5,4.5) {};
    \node[fill=purple!70] at (6.5,4.5) {};
    \node[fill=purple!70] at (0.5,5.5) {};
    \node[fill=purple!70] at (1.5,5.5) {};
    \node[fill=purple!70] at (2.5,5.5) {};
    \node[fill=purple!70] at (3.5,5.5) {};
    \node[fill=purple!70] at (4.5,5.5) {};
    \node[fill=purple!70] at (5.5,5.5) {};
    \node[fill=purple!70] at (0.5,6.5) {};
    \node[fill=purple!70] at (1.5,6.5) {};
    \node[fill=purple!70] at (2.5,6.5) {};
    \node[fill=purple!70] at (3.5,6.5) {};
    \node[fill=purple!70] at (4.5,6.5) {};
    \node[fill=purple!70] at (5.5,6.5) {};
    \node[fill=purple!70] at (6.5,6.5) {};
    \node[fill=purple!70] at (7.5,6.5) {};
    \node[fill=purple!70] at (0.5,7.5) {};
    \node[fill=purple!70] at (1.5,7.5) {};
    \node[fill=purple!70] at (2.5,7.5) {};
    \node[fill=purple!70] at (3.5,7.5) {};
    \node[fill=purple!70] at (4.5,7.5) {};
    \node[fill=purple!70] at (5.5,7.5) {};
    \node[fill=purple!70] at (6.5,7.5) {};
    \node[fill=purple!70] at (7.5,7.5) {};
    \node[fill=purple!70] at (8.5,7.5) {};
\end{tikzpicture}\hspace{0.35cm}\begin{tikzpicture}[scale=0.2,every node/.style={minimum size=0.2cm-\pgflinewidth, outer sep=0pt}]
    \draw[step=1cm,color=black] (0,0) grid (40,8);
    \node[fill=purple!70] at (0.5,0.5) {};
    \node[fill=purple!70] at (1.5,0.5) {};
    \node[fill=purple!70] at (0.5,1.5) {};
    \node[fill=purple!70] at (1.5,1.5) {};
    \node[fill=purple!70] at (2.5,1.5) {};
    \node[fill=purple!70] at (0.5,2.5) {};
    \node[fill=purple!70] at (1.5,2.5) {};
    \node[fill=purple!70] at (2.5,2.5) {};
    \node[fill=purple!70] at (4.5,2.5) {};
    \node[fill=purple!70] at (5.5,2.5) {};
    \node[fill=purple!70] at (6.5,2.5) {};
    \node[fill=purple!70] at (0.5,3.5) {};
    \node[fill=purple!70] at (1.5,3.5) {};
    \node[fill=purple!70] at (2.5,3.5) {};
    \node[fill=purple!70] at (3.5,3.5) {};
    \node[fill=purple!70] at (4.5,3.5) {};
    \node[fill=purple!70] at (5.5,3.5) {};
    \node[fill=purple!70] at (0.5,4.5) {};
    \node[fill=purple!70] at (1.5,4.5) {};
    \node[fill=purple!70] at (2.5,4.5) {};
    \node[fill=purple!70] at (3.5,4.5) {};
    \node[fill=purple!70] at (4.5,4.5) {};
    \node[fill=purple!70] at (5.5,4.5) {};
    \node[fill=purple!70] at (0.5,5.5) {};
    \node[fill=purple!70] at (1.5,5.5) {};
    \node[fill=purple!70] at (2.5,5.5) {};
    \node[fill=purple!70] at (3.5,5.5) {};
    \node[fill=purple!70] at (4.5,5.5) {};
    \node[fill=purple!70] at (0.5,6.5) {};
    \node[fill=purple!70] at (1.5,6.5) {};
    \node[fill=purple!70] at (2.5,6.5) {};
    \node[fill=purple!70] at (3.5,6.5) {};
    \node[fill=purple!70] at (4.5,6.5) {};
    \node[fill=purple!70] at (5.5,6.5) {};
    \node[fill=purple!70] at (0.5,7.5) {};
    \node[fill=purple!70] at (1.5,7.5) {};
    \node[fill=purple!70] at (2.5,7.5) {};
    \node[fill=purple!70] at (3.5,7.5) {};
    \node[fill=purple!70] at (4.5,7.5) {};
    \node[fill=purple!70] at (5.5,7.5) {};
    \node[fill=purple!70] at (6.5,7.5) {};
\end{tikzpicture}
\caption{QoI $Q_1$, $L=3$: models chosen as level $2$ for $TOL_0$ (top left), $TOL_1$ (top right), $TOL_2$ (bottom left) and $TOL_3$ (bottom right).}\label{fig:L3levelsQ1}
\end{figure}
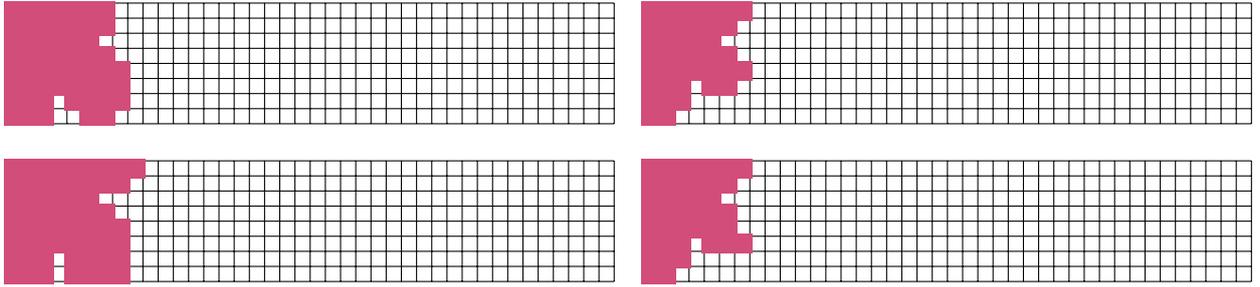

\prev{For $L=2$, the model and level selection procedure chose the model with the finest quadrilateral mesh with elementwise homogenized coefficient ($320$ blocks) as first level, for all tolerances, while the second levels for $TOL_i$, $i=0,\ldots,3$ are shown in Figure \ref{fig:L2levelsQ1} (colored blocks denote the blocks where the microstructure is refined). For $L=3$, the first level corresponds to the finest quadrilateral mesh ($320$ blocks) and the third level to the fine scale model, for all tolerances; the models corresponding to $l=2$ are shown in Figure \ref{fig:L3levelsQ1}. This means that, for $L=2$, the level selection procedure estimated the choice of using the last surrogate model as last level to be cheaper than using a zero bias tolerance and the fine scale model as second level. For $L=3$, instead, using a zero bias tolerance proved to be cheaper than using the bias tolerance associated to the last surrogate selected (we remind that we advance the model selection only until we have a model with tolerance smaller than $\tfrac{TOL_i}{\sqrt 2}$, $i=0,\ldots 3$). The intermediate levels for $L=3$ are different for different tolerances, see Figure \ref{fig:L3levelsQ1}, probably because either choice of these models gives similar costs of the resulting MLMC estimator, and which model is selected depends on the sampled variances. If a coarser or finer model is used at the intermediate level is counterbalanced by taking more or less samples on the finest level, respectively. This can be observed in Table \ref{tab:M_Q3}, which reports, for each tolerance, the number of samples per level used in the MLMC runs (we remind that they do not reflect the total cost of the algorithm, as the cost of the model selection has to be added). If the intermediate model for $L=3$ was always the same, the number of samples on the finest level would roughly quadruple when halving the tolerance; instead, when for one tolerance we have a coarser intermediate model than the previous one (look e.g. $TOL_2$ and $TOL_3$), then the number of samples at $l=3$ is multiplied by a factor larger than $4$, while, when the intermediate model has more refined blocks than in the previous tolerance (e.g. $TOL_1$ compared to $TOL_2$), then the multiplicative factor is lower than $4$. The strong locality of $Q_1$ is reflected by the patterns in the selected models which are highly concentrated around the support of the QoI, and by the fact that, to move from one tolerance to the next one, only few blocks need to be added.}

\begin{table}
\centering
\begin{tabular}{|r|r|r|r|r|r|}
    \hline
    \multirow{2}{*}{Tolerance} &
      \multicolumn{2}{c|}{$L=2$}  &
      \multicolumn{3}{c|}{$L=3$} \\
      \cline{2-6}
    & $\sharp$samples $l=1$ & $\sharp$samples $l=2$ & $\sharp$samples $l=1$ & $\sharp$samples $l=2$ & $\sharp$samples $l=3$\\
    \hline
0.02 & 10163 & 292 & 9176 & 263 & 5\\
0.01 & 46702 & 1253 & 27682 & 916 & 32 \\
0.005 & 128135 & 2905 & 114804 & 2981 & 62 \\
0.0025 & 706639 & 15667 & 566830 & 16787 & 597\\
\hline
  \end{tabular}\caption{\prev{Number of samples per level for $Q_2$.}}\label{tab:M_Q1}
\end{table}

\begin{figure}
\centering
\includegraphics[scale=0.65]{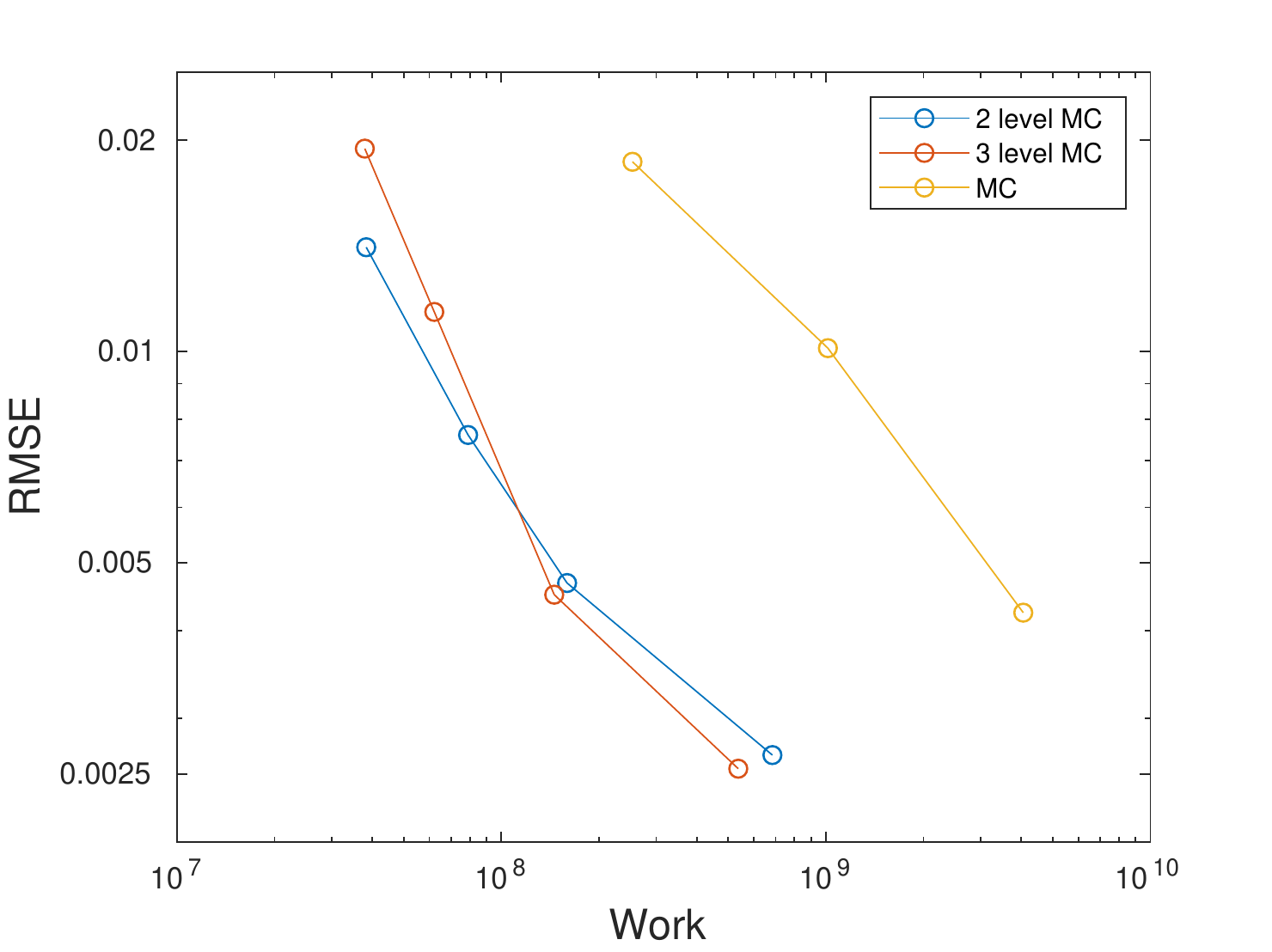}\caption{Convergence plot for the two-level and three level Monte Carlo for $Q_1$, with comparison with the cost of the plain Monte Carlo on the fine scale model.}\label{fig:cvgplotq1}
\end{figure}

\prev{The convergence plot is depicted in Figure \ref{fig:cvgplotq1}. The reference solution is the average of the results of $4$ repetitions of a three-level Monte Carlo with tolerance $0.00125$, where the fine scale model has been used as last level. The root mean squared error has been computed averaging over 30 repetitions. The cost for the standard Monte Carlo algorithm on the fine grid for the first three tolerances is also reported. From Figure \ref{fig:cvgplotq1}, we see that there is no significant gain between the two-level and the three-level Monte Carlo, but there is a gain of a factor of 25 in the cost compared to Monte Carlo.}

\subsection{Plane Strain Elasticity in Random Media}\label{ssec:numexp_ela}

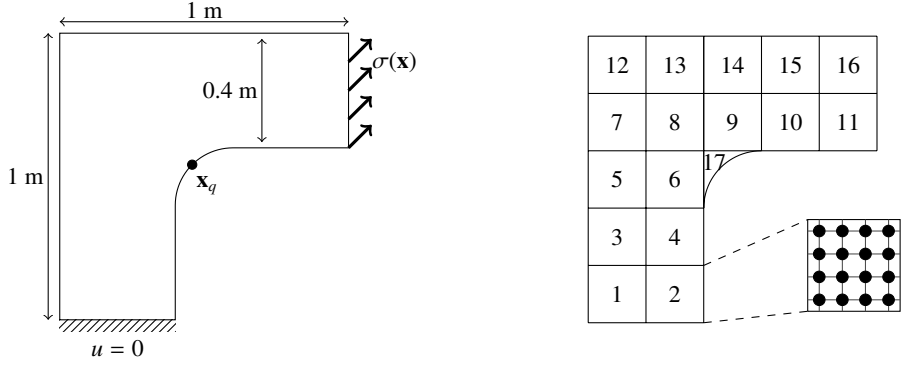
\begin{figure}
 \centering
 \begin{minipage}{0.42\textwidth}
 \centering
\begin{tikzpicture}[scale=0.76]
     \draw [black,domain=90:180] plot ({3+cos(\x)}, {2+sin(\x)}) -- (2,2) -- (2,0)--(0,0)--(0,5)--(5,5)--(5,3)--(3,3);
    \fill[pattern=north east lines] (0,-0.2) rectangle (2,0);
    \node at (1,-0.5) {\footnotesize $u=0$};
   \draw[<->] (-0.2,0)--(-0.2,5) node at (-0.6,2.5) {\footnotesize 1 m};
   \draw[<->] (0,5.2)--(5,5.2) node at (2.5,5.4) {\footnotesize 1 m};
   \draw[<->] (3.5,3.1)--(3.5,4.9) node at (2.95,4) {\footnotesize 0.4 m};
   \draw[->,very thick] (5,3)--(5.4,3.4);
   \draw[->,very thick] (5,3.5)--(5.4,3.9);
   \draw[->,very thick] (5,4)--(5.4,4.4);
   \draw[->,very thick] (5,4.5)--(5.4,4.9);
   \node at (5.8,4.5) {\footnotesize $\mathbf{\sigma(x)}$};
   \filldraw (5*0.4586,5*0.5412) circle (0.08cm);
   \node at (5*0.51,5*0.47) {\footnotesize $\mathbf{x}_q$};
\end{tikzpicture}
 \end{minipage}
 \begin{minipage}{0.42\textwidth}
 \centering
\begin{tikzpicture}[scale=0.76,every node/.style={minimum size=0.7cm-\pgflinewidth, outer sep=0pt}]
    \draw[step=1cm,color=black] (0,0) grid (2,5);
     \draw[step=1cm,color=black] (2,3) grid (5,5);
     \draw [black,domain=90:180] plot ({3+cos(\x)}, {2+sin(\x)}) -- (2,2) -- (2,3);
     \node at (0.5,0.5) {\footnotesize 1};
     \node at (1.5,0.5) {\footnotesize 2};
     \node at (0.5,1.5) {\footnotesize 3};
     \node at (1.5,1.5) {\footnotesize 4};
     \node at (0.5,2.5) {\footnotesize 5};
     \node at (1.5,2.5) {\footnotesize 6};
     \node at (0.5,3.5) {\footnotesize 7};
     \node at (1.5,3.5) {\footnotesize 8};
     \node at (2.5,3.5) {\footnotesize 9};
     \node at (3.5,3.5) {\footnotesize 10};
     \node at (4.5,3.5) {\footnotesize 11};
     \node at (0.5,4.5) {\footnotesize 12};
     \node at (1.5,4.5) {\footnotesize 13};
     \node at (2.5,4.5) {\footnotesize 14};
     \node at (3.5,4.5) {\footnotesize 15};
     \node at (4.5,4.5) {\footnotesize 16};
     \node at (2.18,2.8) {\footnotesize 17};
     \draw (3.8,0.2) rectangle (3.8+1.6,0.2+1.6);
     \draw[white!30!black]  (3.8,0.4) --(3.8+1.6,0.4);
     \draw[white!30!black]  (3.8,0.8) --(3.8+1.6,0.8);
     \draw[white!30!black]  (3.8,1.2) --(3.8+1.6,1.2);
     \draw[white!30!black]  (3.8,1.6) --(3.8+1.6,1.6);
     \draw[white!30!black]  (4,0.2) --(4,1.8);
     \draw[white!30!black]  (4.4,0.2) --(4.4,1.8);
     \draw[white!30!black]  (4.8,0.2) --(4.8,1.8);
     \draw[white!30!black]  (5.2,0.2) --(5.2,1.8);
     \filldraw (4,0.4) circle (0.1cm);
     \filldraw (4,0.8) circle (0.1cm);
     \filldraw (4,1.2) circle (0.1cm);
     \filldraw (4,1.6) circle (0.1cm);
     \filldraw (4.4,0.4) circle (0.1cm);
     \filldraw (4.4,0.8) circle (0.1cm);
     \filldraw (4.4,1.2) circle (0.1cm);
     \filldraw (4.4,1.6) circle (0.1cm);
     \filldraw (4.8,0.4) circle (0.1cm);
     \filldraw (4.8,0.8) circle (0.1cm);
     \filldraw (4.8,1.2) circle (0.1cm);
     \filldraw (4.8,1.6) circle (0.1cm);
     \filldraw (5.2,0.4) circle (0.1cm);
     \filldraw (5.2,0.8) circle (0.1cm);
     \filldraw (5.2,1.2) circle (0.1cm);
     \filldraw (5.2,1.6) circle (0.1cm);
     \draw[dashed] (2,0) -- (3.8,0.2);
     \draw[dashed] (2,1) -- (3.8,1.8);
\end{tikzpicture}
 \end{minipage}
\caption{Setting used for the numerical experiment of Subsection \ref{ssec:numexp_ela}. Left: geometry, boundary conditions and point $P$ around which the QoI is computed. Right: subdivision of the domain into blocks \prev{and zoom of one block to show the subgrid to locate the inclusions.}}\label{fig:Lshapesmooth}
\end{figure}

In this numerical experiment, we test the performance of the error estimator-driven MLMC for \prev{isotropic} plane strain elasticity. The setting is depicted in Figure \ref{fig:Lshapesmooth}. The PDE is \eqref{eq:elastic_forward} with $\mathbf{f}=\mathbf{0}$. We consider homogeneous Dirichlet boundary conditions on the bottom boundary, Neumann boundary conditions with applied traction $\boldsymbol{\sigma}=(500,500)$ at the rightmost boundary and homogeneous Neumann boundary conditions elsewhere; see Figure \ref{fig:Lshapesmooth}.

\prev{We use} a Young's modulus of
\begin{equation}
E(\omega,\bfx) = \begin{cases}
1000 & \text{in the inclusions},\\
100 & \text{in the matrix},
\end{cases}
\end{equation}
(MPa) and a Poisson ratio of $\nu=0.2$ everywhere, as in \cite{odenvemaganti2000}. \prev{These correspond to Lam\'e parameters $\lambda_M = 27.78$(MPa), $\mu_M = 41.67$(MPa) in the matrix and $\lambda_I = 277.78$(MPa), $\mu_I = 416.67$(MPa) in the inclusions.}

\prev{The domain is divided into blocks as in the right plot of Figure \ref{fig:Lshapesmooth}, and the random inclusions are generated as follows:
\begin{itemize}
\item each of the blocks from the first to the sixteenths contains a $n\times n$ tensorial subgrid, denoting the possible positions of the circular inclusions; we use $n=4$;
\item in each square block the number of inclusions, $n_{bl}$, is distributed according to a discrete uniform distribution between some values $0\leq n_{\min}$ and $ n_{\max}\leq n^2$; we use $n_{\min}=0$ and $n_{\max}=n^2$;
\item the average positions of the centers of the inclusions are selected taking the first $n_{bl}$ entries of a random permutation of all possible indices $\left\{1,2,\ldots,n^2\right\}$;
\item the coordinates of the center and the radius of each inclusion are perturbed as in the heat conduction problem, with $h=0.05$(m);
\item for the block $17$, we generate the inclusions as for the other blocks, as if it was a square block; then, we retain the inclusion at the top left corner if this inclusion appears in the square block;
\item all random variables involved are independent.
\end{itemize}}

The QoI is the average trace of the strain around the point $\bfx_q = \left(0.4586,0.5412\right)$ \prev{marked} in the \prev{left} plot of Figure \ref{fig:Lshapesmooth}:
\begin{linenomath}
\begin{equation*}
q_2 =  \frac{1}{\tilde{A}_q}\int_D \chi(\mathbf{x}) \left(\varepsilon_{11}\left(u(\mathbf{x})\right)+\varepsilon_{22}\left(u(\mathbf{x})\right)\right)\,\textit{d}\mathbf{x},
\end{equation*}
\end{linenomath}
with
\begin{linenomath}
\begin{equation*}
\chi(\mathbf{x})=\begin{cases}
1, & \text{for } \lVert  \mathbf{x}-\mathbf{x}_q\rVert\leq r_q, \bfx\in D,\\
\cos^2 \left(\frac{\pi}{2}\frac{\lVert  \mathbf{x}-\mathbf{x}_q\rVert-0.05}{0.05}\right), & \text{for } r_q < \lVert  \mathbf{x}-\mathbf{x}_q\rVert\leq 2r_q, \bfx\in D,\\
0, & \text{for } \lVert  \mathbf{x}-\mathbf{x}_q\rVert> 2r_q, \bfx\in D
\end{cases}\qquad \tilde{A}_q = \int_D \chi(\mathbf{x})\,\textit{d}\mathbf{x},
\end{equation*}
\end{linenomath}
and $r_q=0.05$.

As for the heat conduction problem, we need sufficiently fine discretizations for the modeling error to be the main source of error. We use unstructured, conformal grids with a mesh size of \prev{$0.0045$} in the refined blocks and \prev{$0.045$} in the homogenized region, leading to about \prev{$99391$} degrees of freedom for the \rev{fine scale} model and \prev{$1526$} degrees of freedom for the blockwise homogenized model. 

\prev{Due to Remark \ref{rmk:noscalarg} and since not all blocks have a tensor product structure, as surrogate models we consider the one with blockwise homogenized coefficients and then the ones where the microstructure is resolved in some blocks, as described in Section \ref{sec:heat}, that is $\mathcal{S}=\left\{\mathcal{M}_0\right\}$.}

\prev{We consider the tolerances $TOL_0=0.015, TOL_1=0.0075, TOL_2=0.00375$, for $L=2$ and $L=3$. In Algorithm \ref{alg:modelselect}, we have used $\hat{M}=180$.}

\begin{figure}
\centering
\begin{tikzpicture}[scale=0.63,every node/.style={minimum size=0.62cm-\pgflinewidth, outer sep=0pt}]
    \draw[step=1cm,color=black] (0,0) grid (2,5);
     \draw[step=1cm,color=black] (2,3) grid (5,5);
     \filldraw [purple!70,domain=90:180] plot ({3+cos(\x)}, {2+sin(\x)}) -- (2,2) -- (2,3);
     \node[fill=purple!70] at (1.5,2.5) {};
     \node[fill=purple!70] at (1.5,1.5) {};
     \node[fill=purple!70] at (0.5,2.5) {};
     \node[fill=purple!70] at (1.5,2.5) {};
     \node[fill=purple!70] at (0.5,3.5) {};
     \node[fill=purple!70] at (1.5,3.5) {};
     \node[fill=purple!70] at (2.5,3.5) {};
     \node[fill=purple!70] at (3.5,3.5) {};
     \node[fill=purple!70] at (0.5,4.5) {};
     \node[fill=purple!70] at (1.5,4.5) {};
     \node[fill=purple!70] at (2.5,4.5) {};
     \node[fill=purple!70] at (3.5,4.5) {};
\end{tikzpicture}
\hspace{1cm}
\begin{tikzpicture}[scale=0.63,every node/.style={minimum size=0.62cm-\pgflinewidth, outer sep=0pt}]
    \draw[step=1cm,color=black] (0,0) grid (2,5);
     \draw[step=1cm,color=black] (2,3) grid (5,5);
     \filldraw [purple!70,domain=90:180] plot ({3+cos(\x)}, {2+sin(\x)}) -- (2,2) -- (2,3);
     \node[fill=purple!70] at (1.5,2.5) {};
     \node[fill=purple!70] at (1.5,1.5) {};
     \node[fill=purple!70] at (0.5,1.5) {};
     \node[fill=purple!70] at (0.5,2.5) {};
     \node[fill=purple!70] at (1.5,2.5) {};
     \node[fill=purple!70] at (0.5,3.5) {};
     \node[fill=purple!70] at (1.5,3.5) {};
     \node[fill=purple!70] at (2.5,3.5) {};
     \node[fill=purple!70] at (3.5,3.5) {};
     \node[fill=purple!70] at (0.5,4.5) {};
     \node[fill=purple!70] at (1.5,4.5) {};
     \node[fill=purple!70] at (2.5,4.5) {};
     \node[fill=purple!70] at (3.5,4.5) {};
\end{tikzpicture}
 \caption{QoI $Q_2$, $L=3$: last level for $TOL_0$ (left) and last level for $TOL_1$ and $TOL_2$ (right).}\label{fig:L3levelsQ3}
\end{figure}
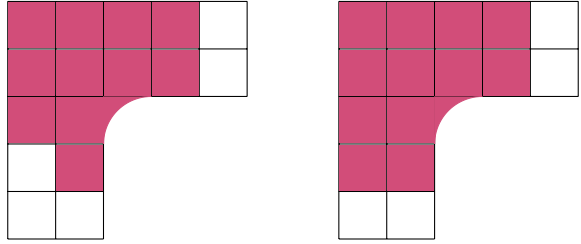

\prev{For $L=2$ and all three tolerances, Algorithm \ref{alg:modelselect} and the level selection chose the two models depicted in Figure \ref{fig:L2levelsQ3}. For $L=3$, level $1$ corresponds to the blockwise homogenized model and level $2$ corresponds to the model on the left of Figure \ref{fig:L2levelsQ3}, for all tolerances; level $3$ is the left model in Figure \ref{fig:L3levelsQ3} for $TOL_0$ and the right model in the same figure for $TOL_1$ and $TOL_2$. Table \ref{tab:M_Q3} shows, for each tolerance, the number of samples per level.}

\begin{figure}[t]
\centering
\includegraphics[scale=0.55]{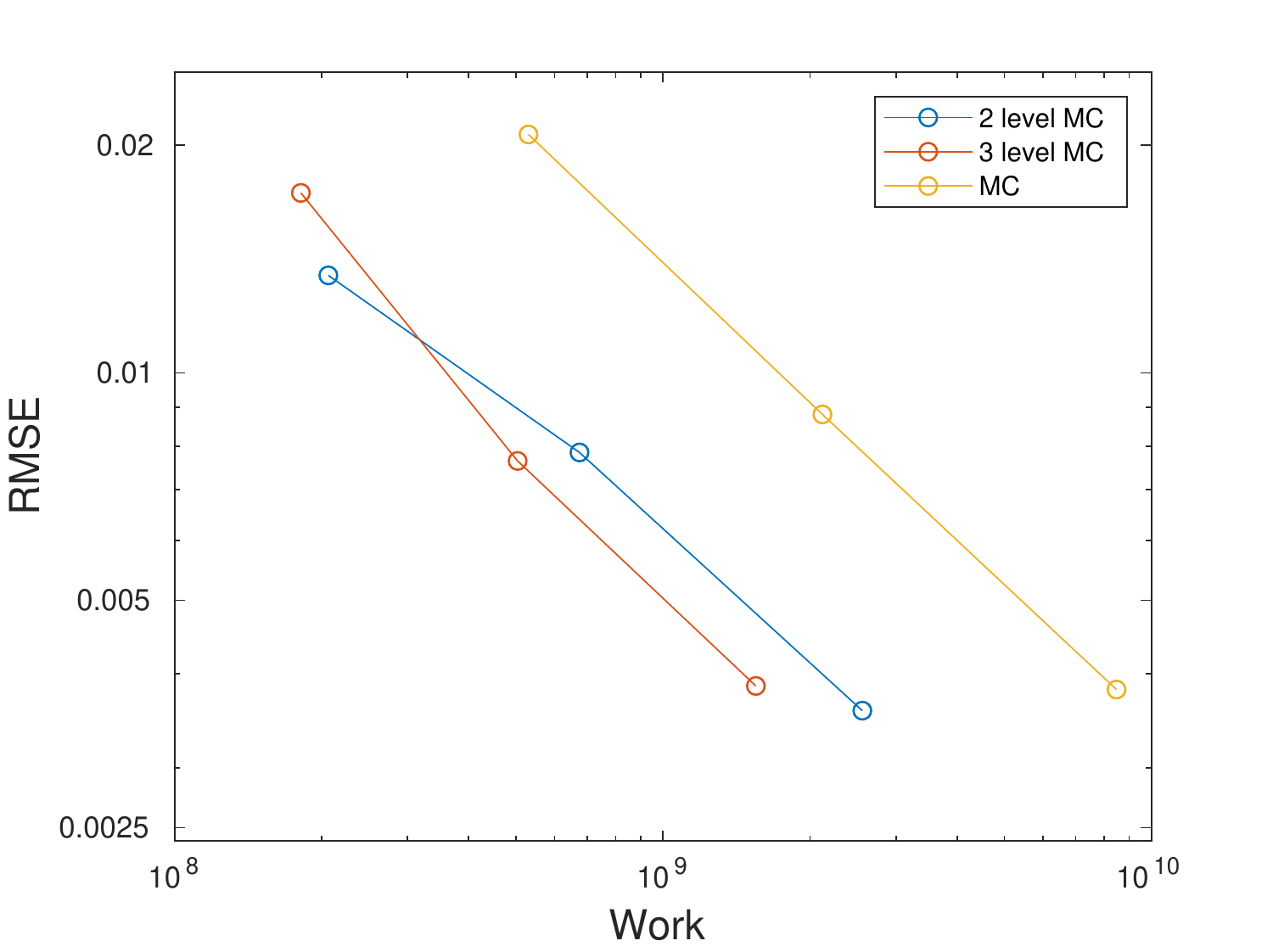}\caption{Convergence plot for the two-level and three level Monte Carlo for $Q_2$, with comparison with the cost of the plain Monte Carlo on the fine scale model.}\label{fig:cvgplotQ3}
\end{figure}

\smallskip
The convergence plot is given in Figure \ref{fig:cvgplotQ3}, \prev{where the cost of the plain Monte Carlo on the fine grid is also reported}. The reference solution has been computed averaging over \prev{$8$} repetitions of a three-level Monte Carlo with tolerance \prev{$0.001875$ (using the fine scale model as last level)} and the root mean squared error has been computed averaging over \prev{25 repetitions}. \prev{The two-level Monte Carlo allows to reduce the cost of plain Monte Carlo by a factor greater than 3, while the three-level Monte Carlo reduces the cost by a factor of almost $6$ and this time it brings savings compared to the two-level case. We explain the more moderate savings of MLMC compared to the previous example with the fact that $Q_2$ is a less local QoI than $Q_1$ (see, for instance, the last levels for $L=3$ in Figure \ref{fig:L3levelsQ3}, where only few blocks are not refined).}

\begin{figure}
\centering
\begin{tikzpicture}[scale=0.62,every node/.style={minimum size=0.62cm-\pgflinewidth, outer sep=0pt}]
    \draw[step=1cm,color=black] (0,0) grid (2,5);
     \draw[step=1cm,color=black] (2,3) grid (5,5);
     \filldraw [purple!70,domain=90:180] plot ({3+cos(\x)}, {2+sin(\x)}) -- (2,2) -- (2,3);
     \node[fill=purple!70] at (1.5,2.5) {};
     \node[fill=purple!70] at (2.5,3.5) {};
\end{tikzpicture}\hspace{1cm}
\begin{tikzpicture}[scale=0.62,every node/.style={minimum size=0.62cm-\pgflinewidth, outer sep=0pt}]
    \draw[step=1cm,color=black] (0,0) grid (2,5);
     \draw[step=1cm,color=black] (2,3) grid (5,5);
     \filldraw [purple!70,domain=90:180] plot ({3+cos(\x)}, {2+sin(\x)}) -- (2,2) -- (2,3);
     \node[fill=purple!70] at (1.5,2.5) {};
     \node[fill=purple!70] at (0.5,0.5) {};
     \node[fill=purple!70] at (1.5,0.5) {};
     \node[fill=purple!70] at (0.5,1.5) {};
     \node[fill=purple!70] at (1.5,1.5) {};
     \node[fill=purple!70] at (0.5,2.5) {};
     \node[fill=purple!70] at (1.5,2.5) {};
     \node[fill=purple!70] at (0.5,3.5) {};
     \node[fill=purple!70] at (1.5,3.5) {};
     \node[fill=purple!70] at (2.5,3.5) {};
     \node[fill=purple!70] at (3.5,3.5) {};
     \node[fill=purple!70] at (4.5,3.5) {};
     \node[fill=purple!70] at (0.5,4.5) {};
     \node[fill=purple!70] at (1.5,4.5) {};
     \node[fill=purple!70] at (2.5,4.5) {};
     \node[fill=purple!70] at (3.5,4.5) {};
     \node[fill=purple!70] at (4.5,4.5) {};
\end{tikzpicture}
 \caption{QoI $Q_2$, levels for $L=2$: level $1$ (left) and level $2$ (right). For all tolerances, the same levels have been \prev{selected}.}\label{fig:L2levelsQ3}     
\end{figure}
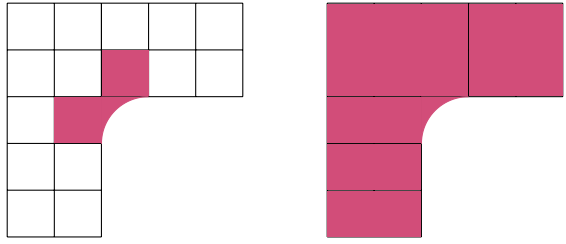

\begin{table}
\centering
\begin{tabular}{|r|r|r|r|r|r|}
    \hline
    \multirow{2}{*}{Tolerance} &
      \multicolumn{2}{c|}{$L=2$}  &
      \multicolumn{3}{c|}{$L=3$} \\
      \cline{2-6}
    & $\sharp$samples $l=1$ & $\sharp$samples $l=2$ & $\sharp$samples $l=1$ & $\sharp$samples $l=2$ & $\sharp$samples $l=3$\\
    \hline
0.015 & 5759 & 241 & 17305 & 2911 & 247\\
0.0075 & 26792 & 912 & 71817 & 11222 & 857\\
0.00375 & 114776 & 3880 & 287374 & 36563 & 3320\\
\hline
  \end{tabular}\caption{\prev{Number of samples per level for $Q_2$.}}\label{tab:M_Q3}
\end{table}


\subsection{Pseudo-3D Linear Elasticity in Random Media}\label{ssec:3d}

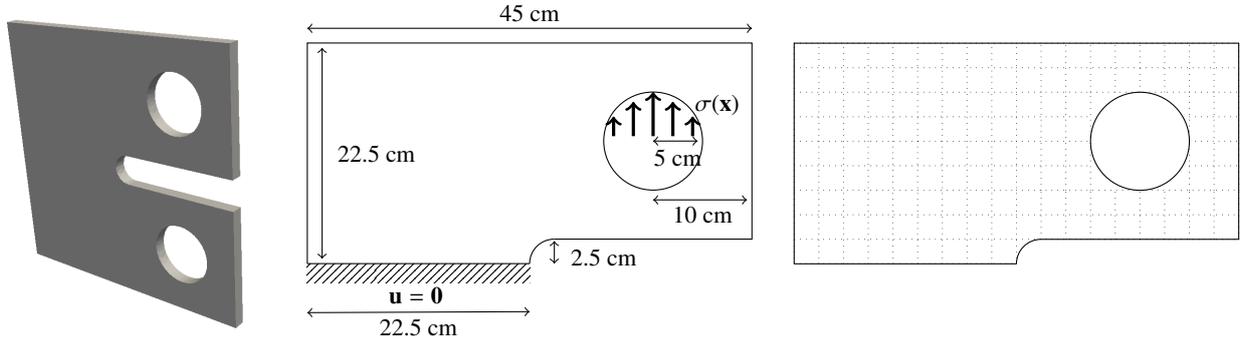
\begin{figure}
\begin{minipage}{0.245\textwidth}
\includegraphics[scale=0.18]{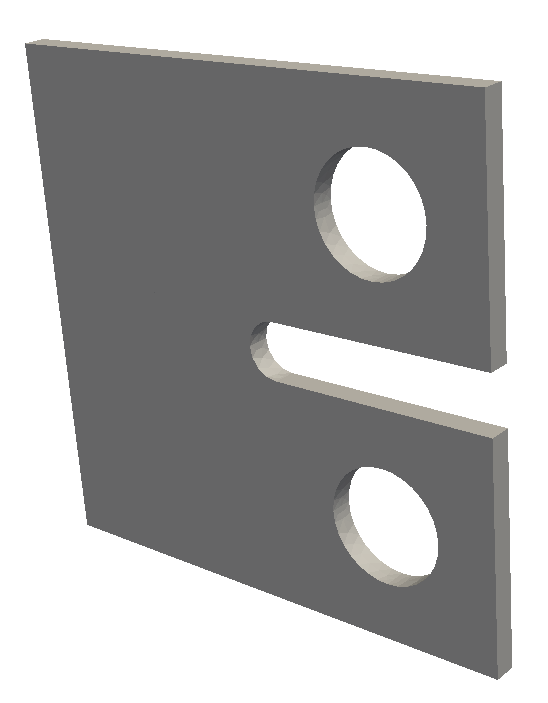}
\end{minipage}
\begin{minipage}{0.385\textwidth}
\begin{tikzpicture}[scale=0.13]
    \draw (0,0)--(0,22.5)--(45,22.5)--(45,2.5)--(25,2.5);
    \draw (0,0)--(22.5,0);
    \draw (35,12.5) circle (5cm);
    \draw (25,2.5) arc (90:180:2.5cm);
    \draw[<->] (0,-5)--(22.5,-4.9) node at (11.25,-6.5) {\footnotesize 22.5 cm};
    \fill[pattern=north east lines] (0,-2) rectangle (22.5,0);
    \node at (11,-3.2) {\footnotesize $\mathbf{u}=\mathbf{0}$};
    \draw[<->] (0,24)--(45,24) node at (22.5,25.5) {\footnotesize 45 cm};
    \draw[<->] (1.5,0.5)--(1.5,22) node at (7,11.25) {\footnotesize 22.5 cm};
    \draw[<->] (35,12.5)--(39.5,12.5) node at (37.5,10.8) {\footnotesize 5 cm};
    \draw[<->] (35,6.5)--(44.5,6.5) node at (40,5) {\footnotesize 10 cm};
    \draw[<->] (25,0)--(25,2.3) node at (30,0.6) {\footnotesize 2.5 cm};
    \draw[->,very thick] (35,13)--(35,17.5);
    \draw[->,very thick] (31,13)--(31,15);
    \draw[->,very thick] (33,13)--(33,16.5);
    \draw[->,very thick] (37,13)--(37,16.5);
    \draw[->,very thick] (39,13)--(39,15);
    \node at (41.5,16.2) {\footnotesize $\mathbf{\sigma(x)}$};
\end{tikzpicture}
\end{minipage}
\begin{minipage}{0.385\textwidth}

\vspace{-0.45cm}
\begin{tikzpicture}[scale=0.13]
    \draw (0,0)--(0,22.5)--(45,22.5)--(45,2.5)--(25,2.5);
    \draw (0,0)--(22.5,0);
    \draw (35,12.5) circle (5cm);
    \draw (25,2.5) arc (90:180:2.5cm);
    \draw[dotted,step=2.5,very thin] (0.0,0.0) grid (22.5,22.5);
    \draw[dotted,step=2.5,very thin] (22.5,2.5) grid (45,7.5);
    \draw[dotted,step=2.5,very thin] (22.5,17.5) grid (45,22.5);
    \draw[dotted,step=2.5,very thin] (22.5,7.5) grid (30,17.5);
    \draw[dotted,step=2.5,very thin] (40,7.5) grid (45,17.5);
    \draw[dotted,very thin] (30,17.5)--(40,17.5);
    \draw[dotted,very thin] (40,7.5)--(40,17.5);
\end{tikzpicture}
\end{minipage}\caption{Setting for the 3D experiment of Subsection \ref{ssec:3d}. Left: compact tension specimen (left), of which we consider the upper half part. Center: geometry in the $(x,y)$-plane and boundary conditions. Right: subdivision into blocks in the $(x,y)$-plane (the four corner blocks around the circular hole have double edge length compared to the other blocks); all blocks have thickness $2$cm in the $z$-direction (as the thickness of the structure).}\label{fig:3dnotched}
\end{figure}

We consider isotropic linear elasticity for the upper half part of the compact tension specimen on the left side of Figure \ref{fig:3dnotched}. The specimen has thickness $2$cm in the $z$-direction, while the geometry in the $(x,y)$-plane is as in the image in the center of Figure \ref{fig:3dnotched}. We consider only one layer of blocks in the $z$-direction, while the right image in Figure \ref{fig:3dnotched} depicts the partition into blocks in the $(x,y)$-plane.

The PDE is \eqref{eq:elastic_forward} with $\mathbf{f}=\mathbf{0}$. We set homogeneous Dirichlet boundary conditions on the bottom boundary, Neumann boundary conditions with applied traction $\boldsymbol{\sigma}=(0,30,0)$ at the upper half surface of the circular hole and homogeneous Neumann boundary conditions elsewhere (see central image in Figure \ref{fig:3dnotched}). Such a setting is common in laboratory tests to determine the crack opening displacement of compact tension specimens.

The inclusions are cylindrical fibers, obtained by generating the cross sections in the $(x,y)$-plane as in the previous subsection, using $n=3$, $n_{\min}=2$ and $n_{\max}=n^2$, and then extending them in the $z$-direction. In the blocks which do not have square cross section, we have proceeded by generating the inclusions as if the cross sections where squares and retaining only five positions in the blocks around the circular hole and one position in the block at the notch; morover, in the blocks around the hole, the center coordinates of the inclusion which is closest to the hole have not been perturbed.

We use a Young's modulus (in (GPa)) and a Poisson ratio of
\begin{linenomath}
\begin{equation}
E(\omega,\bfx) = \begin{cases}
450 & \text{in the inclusions},\\
69 & \text{in the matrix},
\end{cases}\qquad \nu(\omega,\bfx) = \begin{cases}
0.2 & \text{in the inclusions},\\
0.3 & \text{in the matrix}.
\end{cases}
\end{equation}
\end{linenomath}
These correspond to Lam\'e parameters $\lambda_M = 39.81$(GPa), $\mu_M = 26.54$(GPa) in the matrix and $\lambda_I = 125$(GPa), $\mu_I = 187.5$(GPa) in the inclusions.

The QoI is the average dispacement in the $y$-direction (second component of $\mathbf{u}$), denoted by $u_y$, over the surface $S_q:=\left\{45\right\}\times (2.5,22.5)\times (0,2)$:
\begin{linenomath}
\begin{equation}
q_3(\mathbf{u}) = \frac{1}{A_q}\int_{S_q} u_y(\mathbf{x})\,\textit{d}S(\mathbf{x}),
\end{equation}
\end{linenomath}
with $A_q=40$($\text{cm}^2$).

For the spatial discretization, we use a mesh size of $1.25$ in the homogenized region and $0.25$ where the microstructure is resolved, leading to about $15530$ degrees of freedom for the blockwise homogenized model and $556590$ for the fine scale model. As in the previous example, we start the model selection with the blockwise homogenized model, that is $\mathcal{S}=\left\{\mathcal{M}_0\right\}$. The Hashin-Shtrikman bounds are computed as from \eqref{eq:hsela1}, with $d=3$.

We consider the tolerances $TOL_0=0.02$, $TOL_1=0.01$ and $TOL_2=0.007$. For Algorithm \ref{alg:modelselect}, we set $\hat{M}=100$. After the model selection, at each tolerance we have run the level selection procedure for $L=2$ and $L=3$, and then we have chosen the number of levels with minimum cost; for all tolerances, this turned out to be $L=3$. 

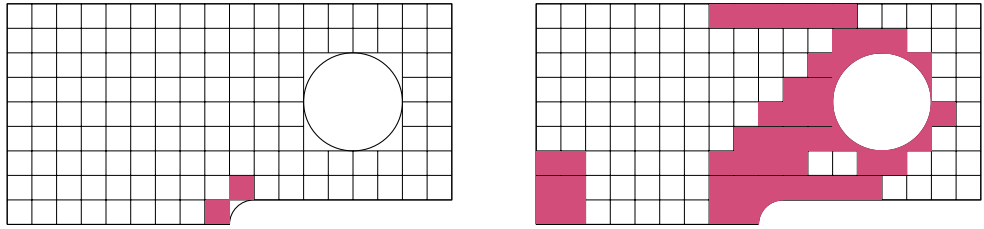
\begin{figure}
\centering
\begin{tikzpicture}[scale=0.13,every node/.style={minimum size=0.33cm-\pgflinewidth, outer sep=0pt}]
    \draw (0,0)--(0,22.5)--(45,22.5)--(45,2.5)--(25,2.5);
    \draw (0,0)--(22.5,0);
    \draw (35,12.5) circle (5cm);
    \draw (25,2.5) arc (90:180:2.5cm);
    \draw[step=2.5,very thin] (0.0,0.0) grid (22.5,22.5);
    \draw[step=2.5,very thin] (22.5,2.5) grid (45,7.5);
    \draw[step=2.5,very thin] (22.5,17.5) grid (45,22.5);
    \draw[step=2.5,very thin] (22.5,7.5) grid (30,17.5);
    \draw[step=2.5,very thin] (40,7.5) grid (45,17.5);
    \draw[very thin] (30,17.5)--(40,17.5);
    \draw[very thin] (40,7.5)--(40,17.5);
    \node[fill=purple!70] at (21.25,1.25) {};
    \node[fill=purple!70] at (23.75,3.75) {};
\end{tikzpicture}\hspace{1cm}
\begin{tikzpicture}[scale=0.13,every node/.style={minimum size=0.33cm-\pgflinewidth, outer sep=0pt}]
    \draw (0,0)--(0,22.5)--(45,22.5)--(45,2.5)--(25,2.5);
    \draw (0,0)--(22.5,0);
    \draw (35,12.5) circle (5cm);coincides with
    \draw (25,2.5) arc (90:180:2.5cm);
    \draw[step=2.5,very thin] (0.0,0.0) grid (22.5,22.5);
    \draw[step=2.5,very thin] (22.5,2.5) grid (45,7.5);
    \draw[step=2.5,very thin] (22.5,17.5) grid (45,22.5);
    \draw[step=2.5,very thin] (22.5,7.5) grid (30,17.5);
    \draw[step=2.5,very thin] (40,7.5) grid (45,17.5);
    \draw[very thin] (30,17.5)--(40,17.5);
    \draw[very thin] (40,7.5)--(40,17.5);
    \node[fill=purple!70] at (1.25,1.25) {};
    \node[fill=purple!70] at (3.75,1.25) {};
    \node[fill=purple!70] at (1.25,3.75) {};
    \node[fill=purple!70] at (3.75,3.75) {};
    \node[fill=purple!70] at (1.25,6.25) {};
    \node[fill=purple!70] at (3.75,6.25) {};
    \node[fill=purple!70] at (21.25,1.25) {};
    \node[fill=purple!70] at (23.75,3.75) {};
    \node[fill=purple!70] at (18.75,1.25) {};
    \node[fill=purple!70] at (18.75,3.75) {};
    \node[fill=purple!70] at (21.25,3.75) {};
    \node[fill=purple!70] at (26.25,3.75) {};
    \node[fill=purple!70] at (28.75,3.75) {};
    \node[fill=purple!70] at (31.25,3.75) {};
    \node[fill=purple!70] at (33.75,3.75) {};
    \node[fill=purple!70] at (18.75,6.25) {};
    \node[fill=purple!70] at (21.25,6.25) {};
    \node[fill=purple!70] at (23.75,6.25) {};
    \node[fill=purple!70] at (26.25,6.25) {};
    \node[fill=purple!70] at (26.25+7.5,6.25) {};
    \node[fill=purple!70] at (36.25,6.25) {};
    \node[fill=purple!70] at (21.25,8.75) {};
    \node[fill=purple!70] at (23.75,8.75) {};
    \node[fill=purple!70] at (26.25,8.75) {};
    \node[fill=purple!70] at (26.25+2.5,8.75) {};
    \node[fill=purple!70] at (23.75,11.25) {};
    \node[fill=purple!70] at (26.25,11.25) {};
    \node[fill=purple!70] at (26.25+2.5,11.25) {};
    \node[fill=purple!70] at (26.25,13.75) {};
    \node[fill=purple!70] at (26.25+2.5,13.75) {};
    \node[fill=purple!70] at (26.25+2.5,16.25) {};
    \node[fill=purple!70] at (26.25+5,18.75) {};
    \node[fill=purple!70] at (26.25+7.5,18.75) {};
    \node[fill=purple!70] at (36.25,18.75) {};
    \node[fill=purple!70] at (26.25+5,18.75+2.5) {};
    \node[fill=purple!70] at (26.25+2.5,18.75+2.5) {};
    \node[fill=purple!70] at (26.25,18.75+2.5) {};
    \node[fill=purple!70] at (21.25+2.5,18.75+2.5) {};
    \node[fill=purple!70] at (21.25,18.75+2.5) {};
    \node[fill=purple!70] at (16.25+2.5,18.75+2.5) {};
    \node[fill=purple!70] at (45-3.75,11.25) {};
    \filldraw [purple!70,domain=90:180] plot ({25+2.5*cos(\x)}, {2.5*sin(\x)}) -- (22.5,0) -- (22.5,2.5);
    \filldraw [purple!70,domain=90:180] plot ({35+5*cos(\x)}, {12.5+5*sin(\x)}) -- (30,12.5) -- (30,17.5);
    \filldraw [purple!70,domain=90:0] plot ({35+5*cos(\x)}, {12.5+5*sin(\x)}) -- (40,17.5) -- (35,17.5);
    \filldraw [purple!70,domain=0:-90] plot ({35+5*cos(\x)}, {12.5+5*sin(\x)}) -- (40,7.5) -- (40,17.5);
    \filldraw [purple!70,domain=-90:-180] plot ({35+5*cos(\x)}, {12.5+5*sin(\x)}) -- (30,7.5) -- (35,7.5);
\end{tikzpicture}\caption{QoI $Q_3$, levels for $L=3$ and $TOL_0$: levels $l=1$ (left) and $l=2$ (center); level $l=3$ coincides with the fine scale model.}\label{fig:3dmod1}
\end{figure}

\begin{figure}
\centering
\begin{tikzpicture}[scale=0.13,every node/.style={minimum size=0.33cm-\pgflinewidth, outer sep=0pt}]
    \draw (0,0)--(0,22.5)--(45,22.5)--(45,2.5)--(25,2.5);
    \draw (0,0)--(22.5,0);
    \draw (35,12.5) circle (5cm);
    \draw (25,2.5) arc (90:180:2.5cm);
    \draw[step=2.5,very thin] (0.0,0.0) grid (22.5,22.5);
    \draw[step=2.5,very thin] (22.5,2.5) grid (45,7.5);
    \draw[step=2.5,very thin] (22.5,17.5) grid (45,22.5);
    \draw[step=2.5,very thin] (22.5,7.5) grid (30,17.5);
    \draw[step=2.5,very thin] (40,7.5) grid (45,17.5);
    \draw[very thin] (30,17.5)--(40,17.5);
    \draw[very thin] (40,7.5)--(40,17.5);
    \node[fill=purple!70] at (1.25,1.25) {};
    \node[fill=purple!70] at (21.25,1.25) {};
    \node[fill=purple!70] at (23.75,3.75) {};
    \node[fill=purple!70] at (18.75,1.25) {};
    \node[fill=purple!70] at (21.25,3.75) {};
    \node[fill=purple!70] at (26.25,3.75) {};
    \node[fill=purple!70] at (23.75,6.25) {};
    \node[fill=purple!70] at (26.25,11.25) {};
    \node[fill=purple!70] at (28.75,13.75) {};
    \node[fill=purple!70] at (23.75,18.75+2.5) {};
    \node[fill=purple!70] at (26.25+2.5,18.75+2.5) {};
    \node[fill=purple!70] at (26.25,18.75+2.5) {};    
    \node[fill=purple!70] at (21.25,18.75+2.5) {};
    \filldraw [purple!70,domain=90:180] plot ({25+2.5*cos(\x)}, {2.5*sin(\x)}) -- (22.5,0) -- (22.5,2.5);
    \filldraw [purple!70,domain=90:180] plot ({35+5*cos(\x)}, {12.5+5*sin(\x)}) -- (30,12.5) -- (30,17.5);
    \filldraw [purple!70,domain=0:-90] plot ({35+5*cos(\x)}, {12.5+5*sin(\x)}) -- (40,7.5) -- (40,17.5);
\end{tikzpicture}\hspace{1cm}
\begin{tikzpicture}[scale=0.13,every node/.style={minimum size=0.33cm-\pgflinewidth, outer sep=0pt}]
    \draw (0,0)--(0,22.5)--(45,22.5)--(45,2.5)--(25,2.5);
    \draw (0,0)--(22.5,0);
    \draw (35,12.5) circle (5cm);
    \draw (25,2.5) arc (90:180:2.5cm);
    \draw[step=2.5,very thin] (0.0,0.0) grid (22.5,22.5);
    \draw[step=2.5,very thin] (22.5,2.5) grid (45,7.5);
    \draw[step=2.5,very thin] (22.5,17.5) grid (45,22.5);
    \draw[step=2.5,very thin] (22.5,7.5) grid (30,17.5);
    \draw[step=2.5,very thin] (40,7.5) grid (45,17.5);
    \draw[very thin] (30,17.5)--(40,17.5);
    \draw[very thin] (40,7.5)--(40,17.5);
    \node[fill=purple!70] at (1.25,1.25) {};
    \node[fill=purple!70] at (21.25,1.25) {};
    \node[fill=purple!70] at (23.75,3.75) {};
    \node[fill=purple!70] at (21.25,3.75) {};
    \node[fill=purple!70] at (26.25,3.75) {};
    \node[fill=purple!70] at (28.75,13.75) {};
    \node[fill=purple!70] at (26.25+2.5,18.75+2.5) {};
    \node[fill=purple!70] at (26.25,18.75+2.5) {};   
    \filldraw [purple!70,domain=90:180] plot ({25+2.5*cos(\x)}, {2.5*sin(\x)}) -- (22.5,0) -- (22.5,2.5);
    \filldraw [purple!70,domain=0:-90] plot ({35+5*cos(\x)}, {12.5+5*sin(\x)}) -- (40,7.5) -- (40,17.5);
\end{tikzpicture}\caption{QoI $Q_3$, level $l=2$ for and $TOL_1$ (left) and $TOL_2$ (right); for both tolerances, level $l=1$ is the model with blockwise homogenized coefficient and level $l=3$ is the fine scale model.}\label{fig:3dmod2}
\end{figure}
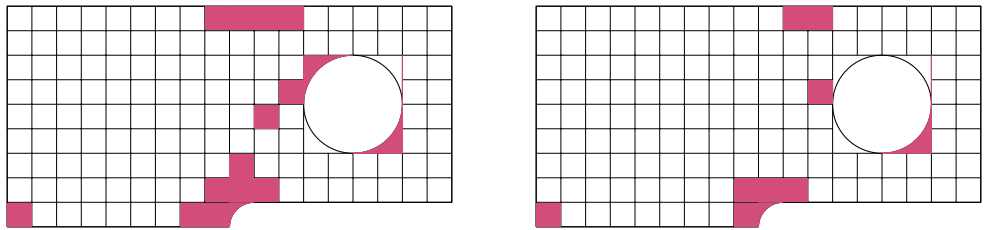

For $L=3$ and $TOL_0$, the first and second levels are shown in Figure \ref{fig:3dmod1}, while the last level is the fine scale model. For $L=3$ and the other two tolerances, the first level is the model with blockwise homogenized coefficient, the last level is the fine scale model and the second level for each of the two tolerances is shown in Figure \ref{fig:3dmod2}. Although the location of the QoI affects the pattern of the blocks where the fine scale is resolved, such local effect is dominated by the local effects at the notch. If we selected the blocks using an a posteriori error estimator for the solution $\mathbf{u}$ itself, probably we would have also started resolving the microstructure around the notch. The number of samples per level for the three-level Monte Carlo runs is reported in Table \ref{tab:M_3d}. 

\begin{figure}[t]
\centering
\includegraphics[scale=0.55]{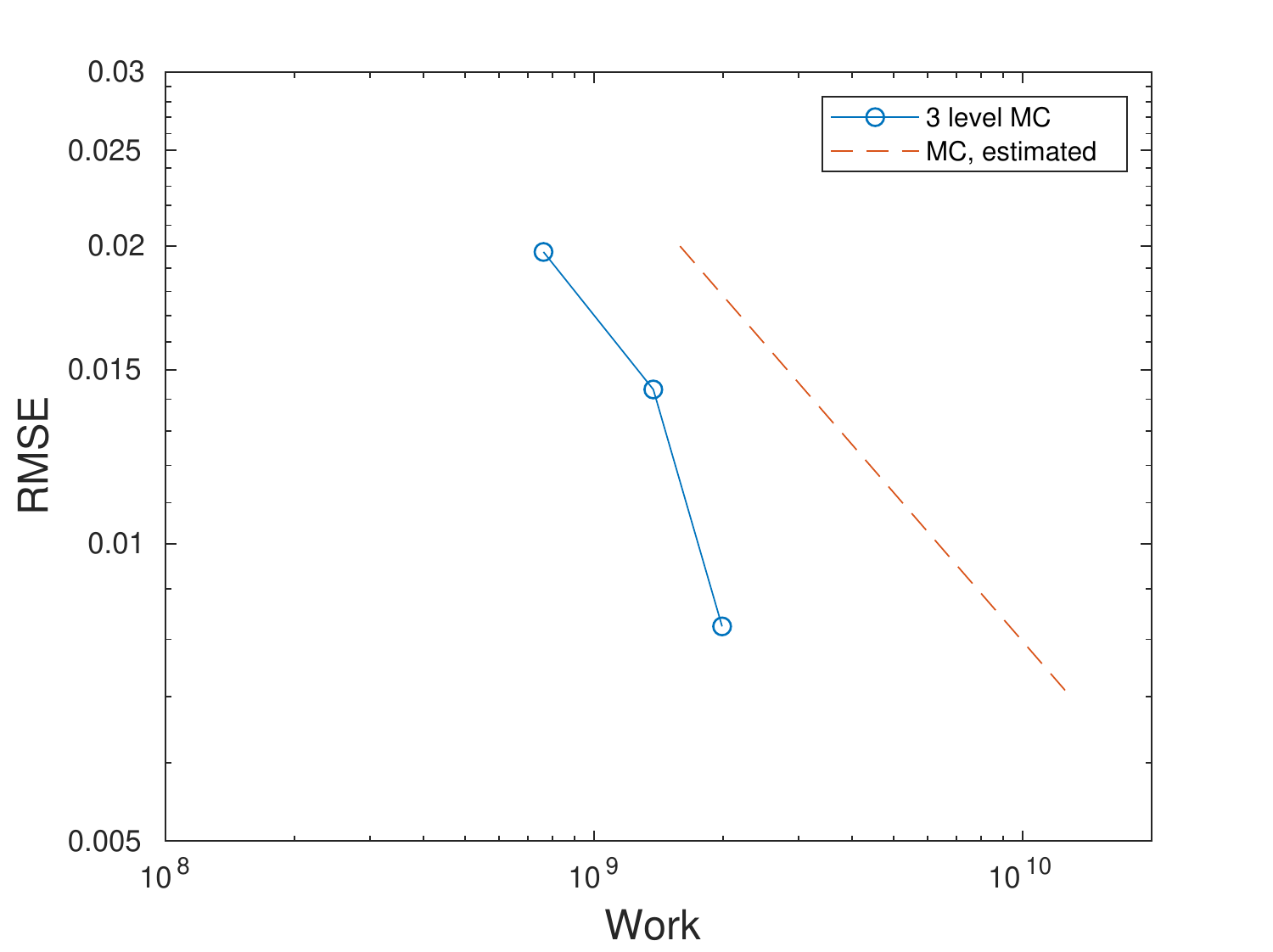}\caption{Convergence plot for three level Monte Carlo for $Q_3$, with comparison with the estimated cost of the plain Monte Carlo on the fine scale model.}\label{fig:3dcvg}
\end{figure}

The convergence plot for the three-level Monte Carlo is shown in Figure \ref{fig:3dcvg}. For each tolerance, the error has been computed from an average of $15$ repetitions, and the reference solution has been computed by averaging the result of $4$ repetitions of the three-level Monte Carlo with tolerance $0.0035$ and the fine scale model as last level. When computing the error for $TOL_2$, we have discarded one of the $15$ repetitions as it was clearly an outlier\footnote{We report the values of the error $|Q^i-Q_{ref}|$ for each repetition $i$ at tolerance $TOL_2$, $Q^i$ denoting the $i$-th repetition and $Q_{ref}$ the reference solution. The errors are as follows, where the outlier is highlighted with bold character: $(0.3406, \mathbf{2.4180}, 0.3214, 0.1394, 0.5990, 1.3066, 1.1595, 0.6635, 0.0963, 0.9966, 0.1403, 1.4047, 0.9330, 0.9980, 0.7539)\cdot 10^{-2}$.}. Since the forward model is very expensive, the cost for plain Monte Carlo is \textsl{estimated}: for each tolerance $TOL_i$, $i=0,1,2$, we have computed the number of samples as $M_i=\tilde{V}/TOL_i^2$, with $\tilde{V}$ a sample average of the variance of the QoI using $500$ samples, and then we have determined an approximation of the cost of Monte Carlo as $W_i = M_i\cdot \tilde{N}$, $i=0,1,2$, with $\tilde{N}$ a sample average of number of degrees of freedom for the fine scale model. From Figure \ref{fig:cvgplotQ3}, we can see that the three-level Monte Carlo allows to have cost savings by a factor of $6$ with respect to plain Monte Carlo. This is for a ratio of $\tfrac{1}{5}$ between the mesh size on the fine scale and the mesh size on the coarse scale. If the ratio was $\frac{1}{10}$ as in the examples of Subsections \ref{ssec:numexp_heat} and \ref{ssec:numexp_ela}, we could expect to gain a factor of around $12$ in the cost.

\begin{table}
\centering
\begin{tabular}{|r|r|r|r|}
    \hline
    \multirow{2}{*}{Tolerance} &
      \multicolumn{3}{c|}{$L=3$} \\
      \cline{2-4}
&$\sharp$samples $l=1$ & $\sharp$samples $l=2$ & $\sharp$samples $l=3$\\
    \hline
0.02 & 5259 & 281 & 58 \\
0.01 & 23601 & 2165 & 414 \\
0.007 & 39327 & 4607 & 843 \\
\hline
  \end{tabular}\caption{Number of samples per level for $Q_3$.}\label{tab:M_3d}
\end{table}

\section{Closing Comments}\label{sec:conclusions}

In this work, a posteriori estimates of modeling error in models of random heterogeneous material are derived and used to construct sequences of surrogate approximations of increasing accuracy of quantities of interest of a \rev{fine scale} base model. This framework provides a basis for new model-based Multilevel Monte Carlo \rev{(}mbMLMC\rev{)} methods. Algorithms for implementing these methods are described and applied to representative examples in stochastic heat conduction and plane elasticity. Numerical experiments indicate that substantial reduction in computational costs can be realized by mbMLMC over standard MC methods. They also indicate that a construction of surrogate models which reflects features of the QoI is essential for a good performance of mbMLMC.

\section*{Acknowledgments}

The work of LS and BW was supported by the German Science Foundation, DFG, WO-671\; 11-1, and by the J. Tinsley Oden Faculty Fellowship Research Program of the ICES Institute (UT Austin).
The work of DF and JTO was supported by the U.S. Department of Energy, Office of Science; Office of Advanced Scientific Computing Research under Award DE-5C0009286.
We benefited from discussion with Robert Lipton on homogenization of random heterogenous media, and we thank Daniel Drzisga for extensive support in the parallelization of the MLMC code.
We also acknowledge with gratitude the early discussions with Kumar Vemaganti of the University of Cincinnati.

\bibliographystyle{abbrv}  
\bibliography{refs}

\begin{thebibliography}{10}

\bibitem{fenicsII}
M.~Aln{\ae}s, J.~Blechta, J.~Hake, A.~Johansson, B.~Kehlet, A.~Logg,
  C.~Richardson, J.~Ring, M.~E. Rognes, and G.~N. Wells.
\newblock {The FEniCS project version 1.5}.
\newblock {\em Archive of Numerical Software}, 3(100):9--23, 2015.

\bibitem{babuskatempone2004}
I.~Babuska, R.~Tempone, and G.~E. Zouraris.
\newblock Galerkin finite element approximations of stochastic elliptic partial
  differential equations.
\newblock {\em SIAM Journal on Numerical Analysis}, 42(2):800--825, 2004.

\bibitem{RBempirical}
M.~Barrault, Y.~Maday, N.~C. Nguyen, and A.~T. Patera.
\newblock An ‘empirical interpolation’ method: application to efficient
  reduced-basis discretization of partial differential equations.
\newblock {\em Comptes Rendus Mathematique}, 339(9):667--672, 2004.

\bibitem{baumanodenprud2009}
P.~T. Bauman, J.~T. Oden, and S.~Prudhomme.
\newblock {Adaptive multiscale modeling of polymeric materials with Arlequin
  coupling and Goals algorithms}.
\newblock {\em Computer Methods in Applied Mechanics and Engineering},
  198(5):799--818, 2009.

\bibitem{bensoussan1978}
A.~Bensoussan, J.-L. Lions, and G.~Papanicolaou.
\newblock {\em Asymptotic analysis for periodic structures}, volume~5.
\newblock North-Holland Publishing Company Amsterdam, 1978.

\bibitem{botevvariance}
Z.~Botev and A.~Ridder.
\newblock Variance reduction.
\newblock {\em Wiley StatsRef: Statistics Reference Online}.

\bibitem{braack2003}
M.~Braack and A.~Ern.
\newblock A posteriori control of modeling errors and discretization errors.
\newblock {\em Multiscale Modeling \& Simulation}, 1(2):221--238, 2003.

\bibitem{buryachenko2007}
V.~Buryachenko.
\newblock {\em Micromechanics of heterogeneous materials}.
\newblock Springer Science \& Business Media, 2007.

\bibitem{Chamoinetal}
L.~Chamoin, J.~T. Oden, and S.~Prudhomme.
\newblock {A stochastic coupling method for atomic-to-continuum Monte-Carlo
  simulations}.
\newblock {\em Computer Methods in Applied Mechanics and Engineering},
  197(43-44):3530--3546, 2008.

\bibitem{cliffe2011}
K.~A. Cliffe, M.~B. Giles, R.~Scheichl, and A.~L. Teckentrup.
\newblock {Multilevel Monte Carlo methods and applications to elliptic PDEs
  with random coefficients}.
\newblock {\em Computing and Visualization in Science}, 14(1):3, 2011.

\bibitem{collier2015contMLMC}
N.~Collier, A.-L. Haji-Ali, F.~Nobile, E.~von Schwerin, and R.~Tempone.
\newblock {A continuation multilevel Monte Carlo algorithm}.
\newblock {\em BIT Numerical Mathematics}, 55(2):399--432, 2015.

\bibitem{eigel2016adaptive}
M.~Eigel, C.~Merdon, and J.~Neumann.
\newblock {An adaptive multilevel Monte Carlo method with stochastic bounds for
  quantities of interest with uncertain data}.
\newblock {\em SIAM/ASA Journal on Uncertainty Quantification},
  4(1):1219--1245, 2016.

\bibitem{giles2008}
M.~B. Giles.
\newblock Multilevel monte carlo path simulation.
\newblock {\em Operations Research}, 56(3):607--617, 2008.

\bibitem{giles_2015}
M.~B. Giles.
\newblock {Multilevel Monte Carlo methods}.
\newblock {\em Acta Numerica}, 24:259–328, 2015.

\bibitem{hoel2014adaptiveMLMC}
H.~Hoel, E.~Von~Schwerin, A.~Szepessy, and R.~Tempone.
\newblock {Implementation and analysis of an adaptive multilevel Monte Carlo
  algorithm}.
\newblock {\em Monte Carlo Methods and Applications}, 20(1):1--41, 2014.

\bibitem{jeulinostoja2001}
D.~Jeulin and M.~Ostoja-Starzewski.
\newblock {\em Mechanics of random and multiscale microstructures}.
\newblock Springer, 2001.

\bibitem{jikovkozlov2012}
V.~V. Jikov, S.~M. Kozlov, and O.~A. Oleinik.
\newblock {\em Homogenization of differential operators and integral
  functionals}.
\newblock Springer Science \& Business Media, 2012.

\bibitem{willcox2017}
H.~Li, V.~Garg, and K.~Willcox.
\newblock {Model Adaptivity for Goal-Oriented Inference using Adjoints}.
\newblock {\em Preprint}, 2017.

\bibitem{fenicsI}
A.~Logg, K.-A. Mardal, and G.~Wells.
\newblock {\em {Automated solution of differential equations by the finite
  element method: The FEniCS book}}, volume~84.
\newblock Springer Science \& Business Media, 2012.

\bibitem{maier2016}
M.~Maier and R.~Rannacher.
\newblock Duality-based adaptivity in finite element discretization of
  heterogeneous multiscale problems.
\newblock {\em Journal of Numerical Mathematics}, 24(3):167--187, 2016.

\bibitem{maier2016arx}
M.~Maier and R.~Rannacher.
\newblock A duality-based optimization approach for model adaptivity in
  heterogeneous multiscale problems.
\newblock {\em arXiv preprint arXiv:1611.09437}, 2016.

\bibitem{mattis2018goal}
S.~A. Mattis and B.~Wohlmuth.
\newblock {Goal-oriented adaptive surrogate construction for stochastic
  inversion}.
\newblock {\em Computer Methods in Applied Mechanics and Engineering},
  339:36--60, 2018.

\bibitem{narayan2014mfsc}
A.~Narayan, C.~Gittelson, and D.~Xiu.
\newblock A stochastic collocation algorithm with multifidelity models.
\newblock {\em SIAM Journal on Scientific Computing}, 36(2):A495--A521, 2014.

\bibitem{odenprudhomme2002}
J.~T. Oden and S.~Prudhomme.
\newblock Estimation of modeling error in computational mechanics.
\newblock {\em Journal of Computational Physics}, 182(2):496--515, 2002.

\bibitem{odenprudrombaum2006}
J.~T. Oden, S.~Prudhomme, A.~Romkes, and P.~T. Bauman.
\newblock {Multiscale modeling of physical phenomena: Adaptive control of
  models}.
\newblock {\em SIAM Journal on Scientific Computing}, 28(6):2359--2389, 2006.

\bibitem{oden1999adaptive}
J.~T. Oden and K.~Vemaganti.
\newblock {Adaptive modeling of composite structures: Modeling error
  estimation}.
\newblock In {\em Texas Institute for Computational and Applied Mathematics}.
  Citeseer, 1999.

\bibitem{odenvemaganti2000}
J.~T. Oden and K.~S. Vemaganti.
\newblock {Estimation of local modeling error and goal-oriented adaptive
  modeling of heterogeneous materials: I. Error estimates and adaptive
  algorithms}.
\newblock {\em Journal of Computational Physics}, 164(1):22--47, 2000.

\bibitem{peherstorfer2016cvg}
B.~Peherstorfer, M.~Gunzburger, and K.~Willcox.
\newblock {Convergence analysis of multifidelity Monte Carlo estimation}.
\newblock {\em Numerische Mathematik}, pages 1--25, 2016.

\bibitem{willcoxgunzburger2016I}
B.~Peherstorfer, K.~Willcox, and M.~Gunzburger.
\newblock {Optimal model management for multifidelity Monte Carlo estimation}.
\newblock {\em SIAM Journal on Scientific Computing}, 38(5):A3163--A3194, 2016.

\bibitem{willcoxgunzburger2016II}
B.~Peherstorfer, K.~Willcox, and M.~Gunzburger.
\newblock Survey of multifidelity methods in uncertainty propagation,
  inference, and optimization.
\newblock {\em SIAM Review}, 60(3):550--591, 2018.

\bibitem{prudhomme2015adaptive}
S.~Prudhomme and C.~M. Bryant.
\newblock {Adaptive surrogate modeling for response surface approximations with
  application to bayesian inference}.
\newblock {\em Advanced Modeling and Simulation in Engineering Sciences},
  2(1):22, 2015.

\bibitem{prudhommeoden1999}
S.~Prudhomme and J.~T. Oden.
\newblock On goal-oriented error estimation for elliptic problems: application
  to the control of pointwise errors.
\newblock {\em Computer Methods in Applied Mechanics and Engineering},
  176(1-4):313--331, 1999.

\bibitem{Romkesetal}
A.~Romkes, J.~T. Oden, and K.~Vemaganti.
\newblock Multi-scale goal-oriented adaptive modeling of random heterogeneous
  materials.
\newblock {\em Mechanics of materials}, 38(8-10):859--872, 2006.

\bibitem{RB}
G.~Rozza, D.~B.~P. Huynh, and A.~T. Patera.
\newblock Reduced basis approximation and a posteriori error estimation for
  affinely parametrized elliptic coercive partial differential equations.
\newblock {\em Archives of Computational Methods in Engineering}, 15(3):1,
  2007.

\bibitem{palencia1983}
E.~Sanchez-Palencia.
\newblock Homogenization method for the study of composite media.
\newblock {\em Asymptotic analysis II}, 985:192--214, 1983.

\bibitem{torquato2013}
S.~Torquato.
\newblock {\em Random heterogeneous materials: microstructure and macroscopic
  properties}, volume~16.
\newblock Springer Science \& Business Media, 2013.

\bibitem{vemagantioden2001}
K.~S. Vemaganti and J.~T. Oden.
\newblock {Estimation of local modeling error and goal-oriented adaptive
  modeling of heterogeneous materials: Part II: a computational environment for
  adaptive modeling of heterogeneous elastic solids}.
\newblock {\em Computer Methods in Applied Mechanics and Engineering},
  190(46):6089--6124, 2001.

\bibitem{vidal2016}
F.~Vidal-Codina, N.~Nguyen, M.~Giles, and J.~Peraire.
\newblock An empirical interpolation and model-variance reduction method for
  computing statistical outputs of parametrized stochastic partial differential
  equations.
\newblock {\em SIAM/ASA Journal on Uncertainty Quantification}, 4(1):244--265,
  2016.

\bibitem{vidal2015}
F.~Vidal-Codina, N.~C. Nguyen, M.~B. Giles, and J.~Peraire.
\newblock {A model and variance reduction method for computing statistical
  outputs of stochastic elliptic partial differential equations}.
\newblock {\em Journal of Computational Physics}, 297:700--720, 2015.

\bibitem{Zaccardi}
C.~Zaccardi, L.~Chamoin, R.~Cottereau, and H.~Ben~Dhia.
\newblock {Error estimation and model adaptation for a stochastic-deterministic
  coupling method based on the Arlequin framework}.
\newblock {\em International Journal for Numerical Methods in Engineering},
  96(2):87--109, 2013.

\bibitem{zhu2017mfmc}
X.~Zhu, E.~M. Linebarger, and D.~Xiu.
\newblock Multi-fidelity stochastic collocation method for computation of
  statistical moments.
\newblock {\em Journal of Computational Physics}, 341:386--396, 2017.

\end{thebibliography}
\index{Bibliography@\emph{Bibliography}}%

%
%
%

\end{document}